# Differential equations over octonions.

Ludkovsky S.V.

15 September 2009

**Abstract**

Differential equations with constant and variable coefficients over octonions are investigated. It is found that different types of differential equations over octonions can be resolved. For this purpose non-commutative line integration is used. Such technique is applied to linear and non-linear partial differential equations in real variables. Possible areas of applications of these results are outlined.

## 1  Introduction.

Differential equations are very important not only for mathematics, but also for their many-sided applications in different sciences including physics, mechanics, other natural sciences, techniques, economics, etc. [2, 7, 13, 34, 16, 27, 39]. Mainly differential equations are considered over fields such as real, complex, or with non-archimedean norms. Recently they are also begun to be studied over Clifford algebras [9, 10, 11].

This article is devoted to differential equations of octonion variables for octonion valued functions. It opens new possibilities for finding solutions or integrating differential equations not only ordinary of real of complex variables, but also some types of partial differential equations including non-linear and their systems. For this purpose previous results of the author on non-commutative line integration over Cayley-Dickson algebras are used [22, 21, 23]. In these works the super-differentiable functions of the Cayley-Dickson variable $z$ were investigated as well. Such functions were also called $z$-holomorphic or $z$-differentiable, or shortly holomorphic or differentiable, when the Cayley-Dickson algebra is outlined.

As it is well-known, quaternions were first introduced by W.R. Hamilton in the second half of the 19-th century. Hamilton had planned to use them for problems of mechanics and mathematics [12, 33]. Later on Cayley had developed their generalization known as the octonion algebra. Then Dickson had investigated their generalizations known now as the Cayley-Dickson algebras [1, 4, 14].

Octonions and quaternions are widely used in physics, but mainly algebraically. Maxwell had utilized quaternions to derive his equations of electrodynamics, but then he had rewritten them in real coordinates.

Already Yang and Mills had used them in quantum field theory, but theory of functions over octonions and quaternions in their times was not sufficiently developed to satisfy their needs. They have formulated the problem of developing analysis over octonions and quaternions [8]. This is natural, because quantum fields are frequently non-abelian [36].

Dirac had used complexified bi-quaternions to solve spin problems of quantum mechanics. He had solved the Klein-Gordon hyperbolic partial differential equation with constant coefficients. The technique of this article developed below permits to integrate also partial differential equations of elliptic and hyperbolic types with variable coefficients. The latter equations are important in many applications, for example, in physics, mechanics, quantum mechanics, in economics for describing processes in the exchange stock market as well.

Each Cayley-Dickson algebra $\mathcal{A}_{r+1}$ is obtained from the preceding $\mathcal{A}_r$ with the help of the so called doubling procedure [1, 14, 17]. This gives the family of embedded algebras: $\mathcal{A}_r \hookrightarrow \mathcal{A}_{r+1} \hookrightarrow \dots$. For a unification of notation it is convenient to put: $\mathcal{A}_0 = \mathbf{R}$ for the real field, $\mathcal{A}_1 = \mathbf{C}$ for the complex field, $\mathcal{A}_2 = \mathbf{H}$ denotes the quaternion skew field, $\mathcal{A}_3 = \mathbf{O}$ is the octonion algebra, $\mathcal{A}_4$ denotes the sedenion algebra. The quaternion skew field is associative, but non-commutative. The octonion algebra is the alternative division algebra with the multiplicative norm. The sedenion algebra and Cayley-Dickson algebras of higher order $r \geq 4$ are not division algebras and have not any non-trivial multiplicative norm. Each equation of the form $ax = b$ with non-zero octonion $a$ and any octonion $b$ can be resolved in the octonion algebra: $x = a^{-1}b$, but it may be non-resolvable in Cayley-Dickson algebras of higher order $r \geq 4$ because of divisors of zero. Therefore, in this article differential equations are considered with octonion or quaternion variables for octonion or quaternion valued functions. If a given differential equation permits a generalization for Cayley-Dickson algebras of higher order it will be indicated.

We recall the doubling procedure for the Cayley-Dickson algebra $\mathcal{A}_{r+1}$ from $\mathcal{A}_r$, because it is frequently used. Each Cayley-Dickson number $z \in \mathcal{A}_{r+1}$ is written in the form $z = \xi + \eta\mathbf{l}$, where $\mathbf{l}^2 = -1$, $\mathbf{l} \notin \mathcal{A}_r$, $\xi, \eta \in \mathcal{A}_r$. The addition of such numbers is componentwise. The conjugate of any Cayley-Dickson number $z$ is given by the formula:

(1) $z^* := \xi^* - \eta\mathbf{l}$.

The multiplication in $\mathcal{A}_{r+1}$ is prescribed by the following equation:

(2) $(\xi + \eta\mathbf{l})(\gamma + \delta\mathbf{l}) = (\xi\gamma - \tilde{\delta}\eta) + (\delta\xi + \eta\tilde{\gamma})\mathbf{l}$

for each $\xi$, $\eta$, $\gamma$, $\delta \in \mathcal{A}_r$, $z := \xi + \eta\mathbf{l} \in \mathcal{A}_{r+1}$, $\zeta := \gamma + \delta\mathbf{l} \in \mathcal{A}_{r+1}$.

The basis of $\mathcal{A}_r$ over $\mathbf{R}$ is denoted by $\mathbf{b}_r := \mathbf{b} := \{1, i_1, ..., i_{2^r-1}\}$, where $i_s^2 = -1$ for each $1 \leq s \leq 2^r - 1$, $i_{2^r} := \mathbf{l}$ is the additional element of the doubling procedure of $\mathcal{A}_{r+1}$ from $\mathcal{A}_r$, choose $i_{2^r+m} = i_m\mathbf{l}$ for each $m = 1, ..., 2^r - 1$, $i_0 := 1$. This implies that $\xi\mathbf{l} = \mathbf{l}\xi^*$ for each $\xi \in \mathcal{A}_r$, when $1 \leq r$.

The technique developed below permits to integrate ordinary differential equations of different types and their systems of ordinary differential equations up to 8 equations. An octonion variable is 8 dimensional over the real field, so (super-)differential equations over octonions encompass some partial differential equations of up to 8 real variables. This gives effective ways for their integration with the help of a non-commutative non-associative line integral over octonions. The quaternion skew field is the subalgebra of the octonion algebra and the line integral over quaternions is already associative.

Section 2 of this article contains the theorem about an existence of solutions of ordinary and partial differential equations over octonions.

Differential equations over octonions treated below in Section 3 include **R**-linear octonion additive differential operators, some types of non-linear equations including that of Bernoulli. Moreover, differential equations in complete differentials and their integration multipliers are studied. Differential equations of $n$-th order are investigated as well in Section 3. Examples illustrating this new technique for ordinary and some types of partial differential equations are given. It is shown, that the new technique over octonions enlarges possibilities of the integration of ordinary and partial linear or non-linear differential equations in comparison with the real and complex fields.

Main results of this paper are obtained for the first time.

## 2  Existence of solutions.

**1. Remark.** Suppose there is given a system of partial differential equations:

$$(1) \quad \partial^{n_j} u_j/\partial t^{n_j} = F_j(t, x_1, ..., x_n, u_1, ..., u_N, ..., \partial^k u_l/\partial t^{k_0} \partial x_1^{k_1}...\partial x_n^{k_n}, ...),$$

where $j, l = 1, ..., N$; $k_0 + k_1 + ... + k_n = k \leq n_l$; $k_0 < n_l$, $N, n \in \mathbf{N}$, $0 \leq k \in \mathbf{Z}$, $(x_1, ..., x_n) \in U$, $x_j \in \mathbf{K}$ for each $j$, $U$ is a domain in $\mathbf{K^n}$, where $\mathbf{K} = \mathbf{H}$ is either the quaternion skew field or $\mathbf{K} = \mathbf{O}$ is the octonion algebra. When a number $t$ is real the $F_j$ are $\mathbf{K}$-valued functions or functionals. In the case $t \in \mathbf{K}$ each $F_j$ is a mapping with values in the space $L_q(\mathbf{K}^{\otimes n_j}; \mathbf{K})$ of all $\mathbf{R}$ $n_j$-poly-homogeneous $\mathbf{K}$-additive operators with values in $\mathbf{K}$. Each partial derivative $\partial^k u_l/\partial t^{k_0} \partial x_1^{k_1}...\partial x_n^{k_n}$ is the mapping with values in $L_q(\mathbf{K}^{\otimes s}; \mathbf{K})$, where $s = k - k_0$ for $t \in \mathbf{R}$ and $s = k$ for $t \in \mathbf{K}$, $\mathbf{K} = \mathbf{H}$ or $\mathbf{K} = \mathbf{O}$.

For $t \in \mathbf{R}$ there are given initial conditions:

$$(2) \quad (\partial^s u_j/\partial t^s)|_{t=t_0} = \phi_j^{(s)}(x_1, x_2, ..., x_n)$$

while for $t \in \mathbf{K}$ let

$$(2') \quad (\partial^s u_j/\partial t^s)|_{Re(t)=Re(t_0)} = \phi_j^{(s)}(t; x_1, x_2, ..., x_n)$$

for each $s = 0, 1, ..., n_j - 1$ and each $(x_1, ..., x_n) \in U$, where $\partial^0 u_j/\partial t^0 := u_j$, in addition each $\phi_j^{(s)}(t; x_1, x_2, ..., x_n)$ is locally analytic by all variables and independent from $Re(t)$ (see also the introduction).

We search a solution of the problem $(1, 2)$ in a neighborhood $V$ of the given point $(t_0, x_1^0, ..., x_n^0)$ either in $\mathbf{R} \times U$ or in $\mathbf{K} \times U$. For our convenience we introduce the notation:

$$(\partial^{k-k_0}\phi_j^{(k_0)}/\partial x_1^{k_1}...\partial x_n^{k_n})|_{(x_1=x_1^0,...,x_n=x_n^0)} = \phi_{(j,k_0,k_1,...,k_n)}^0$$

for each $j = 1, ..., N$, $k_0 + k_1 + ... + k_n = k \leq n_j$.

(3). Suppose that all functions $F_j$ are analytic by all variables and holomorphic by $x_1, ..., x_n$; $u_1, ..., u_N$, also holomorphic by $t$, when $t \in \mathbf{K}$, in a neighborhood $Y$ of the point $(t_0, x_1^0, ..., x_n^0, ..., \phi_{j,k_0,k_1,...,k_n}^0, ...)$, where analyticity by operators $A_{l,k} := \partial^k u_l/\partial t^{k_0}\partial x_1^{k_1}...\partial x_n^{k_n}$ is defined relative to their suitable compositions and the convergence of a series of $F_j$ by $\{A_{l,k} : l, k\}$ is supposed to be relative to the operator norm, where $A_{l,k} \in \bigcup_{m=1}^{\infty} L_q(\mathbf{K}^{\otimes m}, \mathbf{K})$. Let also all functions $\phi_j^{(k)}$ be $\mathbf{K}$-holomorphic (i.e. locally $x$-analytic) in a neighborhood $X$ of $(x_1^0, ..., x_n^0)$.

**2. Theorem.** *The problem $1(1-3)$ has a holomorphic solution either by all variables $(t, x_1, ..., x_n)$ in some neighborhood of $(t_0, x_1^0, ..., x_n^0)$, when $t \in \mathbf{K}$, or holomorphic by the variables $(x_1, ..., x_n)$ and (locally) analytic by $t$, when $t \in \mathbf{R}$. This solution is unique in such class of functions.*

**Proof.** To prove this theorem we introduce new functions, when $n_j > 1$, so that $v_{j,p} := \partial^p u_j/\partial t^p$ for each $p = 1, ..., n_j - 1$, $w_{l,k} := \partial^k u_l/\partial t^{k_0}\partial x_1^{k_1}...\partial x_n^{k_n}$ for each $k_0 + k_1 + ... + k_n = k \leq n_l$ and $k_0 < n_l$. Then $\partial v_{j,p}/\partial t = v_{j,p+1}$ for each $p < n_j$ and $\partial v_{j,n_j-1}/\partial t = F_j$, also $w_{j,k} = \partial^k v_{l,k_0}/\partial x_1^{k_1}...\partial x_n^{k_n}$ for each $(k_0, ..., k_n)$ and $l$, where $F_j = F_j(t, x_1, ..., x_n, u_1, ..., u_N, ..., w_{j,k}, ...)$. The initial conditions now become:

$$u_j|_{Re(t)=Re(t_0)} = \phi_j^{(0)}(t; x_1, x_2, ..., x_n),$$

$$v_{j,s}|_{Re(t)=Re(t_0)} = \phi_j^{(s)}(t; x_1, x_2, ..., x_n)$$

for each $s = 1, ..., n_j - 1$ and each $(x_1, ..., x_n) \in U$ and every $t \in \mathbf{K}$ with $Re(t) = Re(t_0)$ in the considered domain. We can take as a neighborhood of $t_0$ a ball in $\mathbf{K}$ with center at $t_0$ of a sufficiently small positive radius. In the real case of the variable $t$ it would be an interval in $\mathbf{R}$. For $t \in \mathbf{R}$ we automatically get initial conditions from these general.

If $\partial(v_{j,k-k_0} - \partial^{k-k_0}u/\partial x_1^{k_1}...\partial x_n^{k_n})/\partial t = 0$ for each $t$ and $(x_1, ..., x_n) \in U$, then $(v_{j,k-k_0} - \partial^{k-k_0}u/\partial x_1^{k_1}...\partial x_n^{k_n})$ is independent from $t$. If in addition, $v_{j,k-k_0}|_{Re(t)=Re(t_0)} = \phi_j^{(k-k_0)}(t; x_1, x_2, ..., x_n)$ for each $(x_1, ..., x_n) \in U$ and each $t$ for the given domain, then $\partial^{k-k_0} u/\partial x_1^{k_1}...\partial x_n^{k_n} = \phi_j^{(k-k_0)}(t; x_1, x_2, ..., x_n)$ in $Y$ due to analyticity of functions.

Therefore, this procedure of the introduction of new functions $v_{j,s}$ reduces the problem $1(1-3)$ to the problem, when $n_j = 1$ for each $j$. It is possible to simplify further the problem, making substitutions: $w_j(t,x) := u_j(t,x) - \phi_j^{(0)}(t;x)$, $w_{j,s}(t,x) := v_{j,s}(t,x) - \phi_j^{(s)}(t;x)$ for each $j, s$, where $x = (x_1, ..., x_n)$. For new functions $w_j$ and $w_{j,s}$ the problem has the same form, but with zero initial conditions. Thus without loss of generality consider zero initial conditions and denote $w_j$ and $w_{j,s}$ again by $u_j$ and $v_{j,s}$.

Consider at first the case when each $F_j$ is **R**-homogeneous and **K**-additive by each $u_j$ and each $v_{j,s}$, where $n_j = 1$ for each $j$. Without loss of generality consider $t \in \mathbf{K}$, since when $t \in \mathbf{R}$, then due to the local analyticity by $t$ it is possible to take the extension of each $F_j$ to $F_j \in L_q$ on a neighborhood of $t_0$ in **K**. In view of Proposition 2.18 and Remark 2.18.1 in [21, 22] and Theorem 3.2.14 [25] if a solution of the problem exists in the considered class of holomorphic functions $u_j$, $j = 1, ..., N$, then it is unique.

Consider now series of $F_j$, $u_j$, $v_{j,s}$ for each $j$ and $s$ by all their variables $(t - t_0, x_1 - x_1^0, ..., x_n - x_n^0)$ in a neighborhood of $(t_0, x_1^0, ..., x_n^0)$. Expansion coefficients of $u_j$ and $v_{j,s}$ can be evaluated from the system of equations by induction in the order of increasing powers, since $u_j(t,x)|_{Re(t)=Re(t_0)} = 0$ and $v_{j,s}(t,x)|_{Re(t)=Re(t_0)} = 0$ for each $j, s$ and each $x = (x_1, ..., x_n) \in U$ and every $t$ with the prescribed real part $Re(t) = Re(t_0)$ in the considered domain. Now prove the local convergence of the series.

It is possible to make the shift of variables such that one can consider $(t_0, x_1^0, ..., x_n^0) = 0$ without loss of generality. Let the series decomposition $g(t,x) = \sum_{|k|\geq 0}\{(a_k, (t,x)^k)\}_{q(2|k|)}$ be converging at each point $b = (b_0, b_1, ..., b_n)$ with given values $|b_0|,...,|b_n|$, where $b_l \neq 0$ for each $l = 0, 1, ..., n$. Then there exists a constant $S > 0$ such that $|\sum_{k=(k_0,...,k_n)}\{(a_k, (t,x)^k)\}_{q(2|k|)}| \leq S$ for each $k$, $|k| := k_0 + ... + k_n$. Therefore, this implies the inequality:

$$|\sum_{k=(k_0,...,k_n)} \{(a_k, (\beta_0, \beta_1, ..., \beta_n)^k)\}_{q(2|k|)}| \leq S|b_0|^{-k_0}...|b_n|^{-k_n}$$

for each $\beta_0 \in \mathbf{K}, ..., \beta_n \in \mathbf{K}$ with $|\beta_0| = 1,...,|\beta_n| = 1$ and each integers $0 \leq k_0, ..., k_n < \infty$, where vectors of associators $q(*)$ indicating on orders of brackets or multiplications are important only over $\mathbf{O}$, $a_k = a_k(g)$ (see Formulas $2.13(ii-iv)$ in [25]). Take the function $Q(t,x) := S(1-t/|a_0|)^{-1}...(1-x_n/|a_n|)^{-1}$. Then $Q(t,x)$ is characterized by the condition:

$$|\sum_{k=(k_0,...,k_n)} \{(a_k(g),(\beta_0,\beta_1,...,\beta_n)^k)\}_{q(2|k|)}| \leq |\sum_{k=(k_0,...,k_n)} \{(a_k(Q),(\beta_0,\beta_1,...,\beta_n)^k)\}_{q(2|k|)}|$$

for each $k$ and each $\beta_0,...,\beta_n$ in the algebra $\mathbf{K}$ with the natural restrictions $|\beta_0|=1,...,|\beta_n|=1$. The multiplications $z \mapsto \beta_j^p z$ and $z \mapsto z\beta_j^p$ are pseudoconformal mappings from $\mathbf{K}$ onto $\mathbf{K}$ for each $p \in \mathbf{N}$. Such holomorphic function $Q$ is called the majorant of the holomorphic function $g$. In the case of variables $A_{j,s} \in L_q$ consider majorants in the class of analytic functions by $A_{j,s}$ and their suitable compositions, where we choose $b_{j,s} \in \mathbf{K}$ such that $|b_{j,s}| := \|A_{j,s}\|$, since we can take $A_{j,s}.(h_1,...,h_s)$ and $F_j.h_0$ for $F_j \in L_1(\mathbf{K},\mathbf{K})$ for arbitrary $h_0, h_1,...,h_s \in \mathbf{K}$ of unit norms $|h_0|=1$, $|h_1|=1,...,|h_s|=1$. Then such composite functions while action on $h_0, h_1,..., h_s$ have also series expansions which are uniformly convergent series by $h_0, h_1,...,h_s \in B(\mathbf{K},0,1)$.

Let the problem I be the initial problem formulated in §1 and let the problem II be the problem, in which each participating function is substituted by its majorant. Let $U_j$, $V_{j,s}$ denotes the solution of problem II. Then

$$|\sum_{k=(0,k_1,...,k_n)} \{(a_k(v_{j,s}),(\beta_0,\beta_1,...,\beta_n)^k)\}_{q(2|k|)}| \leq |\sum_{k=(0,k_1,...,k_n)} \{(a_k(V_{j,s}),(\beta_0,\beta_1,...,\beta_n)^k)\}_{q(2|k|)}|$$

for each $\beta_0,...,\beta_n$ in $\mathbf{K}$ with $|\beta_0|=1,...,|\beta_n|=1$ and each $j$ and $s$ due to the initial conditions, where $V_{j,0} := U_j$ and $v_{j,0} := u_j$. In the case $k_0 > 0$ coefficients $a_k(v_{j,s})$ or $a_k(V_{j,s})$ are obtained from coefficients $a_m(v_{j,s})$ or $a_m(V_{j,s})$ respectively with the help of addition and multiplication with numbers subjected to the conditions $0 \leq m_0 < k_0$. Therefore, by induction these inequalities are true for each $k$. Thus if problem II has a solution, the problem I also has a solution.

Choose constants $S > 0$ and $a > 0$ such that the function $S[1-(t/\alpha + x_1+...+x_n)/a]^{-1}$ for some $0 < \alpha < 1$ would be the majorant for all terms of the system besides free terms, as well as take $T[1-(t/\alpha+x_1+...+x_n)/a]^{-1}$ as majorant for free terms, since such majorants with $S = S_j$, $a = a_j$, $T = T_j$ can be taken for each $F_j$. Hence we can take $S = \max_{j=1,...,n} S_j$, $T := \max_{j=1,\dots,n} T_j$, $a = \min_{j=1,\dots,n} a_j$, where $\alpha$ we shall choose below. For it problem II can be taken as:

$\partial U_j / \partial t = [\sum_{j,k} \partial U_j / \partial x_k + \sum_j U_j I + CI] S[1 - (t/\alpha + x_1 + \dots + x_n)/a]^{-1}$,

where $C = T/S$, while $I \in L_1$ is the constant unit operator. Seek a solution $U_j(t,x) = U(z)$, where $z = t/\alpha + x_1 + \dots + x_n$. Then we get the equation:

$\alpha^{-1} dU(z)/dz = (NndU(z)/dz + NUI + CI)S(1 - z/a)^{-1}$.

Choose a number $\alpha > 0$ sufficiently large such that $1/\alpha - NnS(1-z/a)^{-1} \neq 0$ in a neighborhood of $z = 0$. Then the problem is:

$(\alpha^{-1} - S(1-z/a)^{-1} Nn) dU(z)/dz = I(NU + C)$, hence

$(NU+C)^{-1}(dU(z)/dz) = Ig(z)$, where $g(z) = (\alpha^{-1} - NnS(1-z/a)^{-1})^{-1}$

due to the alternativity of **K**. The integration of the latter equation gives the $z$-differentiable (holomorphic) solution due to Theorem 3.2.14 [25] and Remark 2.18.1 [21, 22], since the mapping $g$ is $z$-differentiable. Now estimate its expansion coefficients. For this consider an algebraic embedding $\theta$ of the complex field **C** into **K**. The restriction of $U(z)$ on this copy of **C** has the form:

$U(z) = [\exp((\int_0^z y(\xi) d\xi) N/C) - 1] C/N$, for each $z \in \theta(\mathbf{C})$,

where $y(z) = (\alpha^{-1} - NnS(1-z/a)^{-1})^{-1} CS(1-z/a)^{-1}$. The function $S(1-z/a)^{-1}$ has the nonnegative expansion coefficients by powers $z^k$. Hence $y(\xi)$ also has the nonnegative expansion coefficients by $\xi$. Thus the indefinite integral

$$s(z) = N \int_0^z y(\xi) d\xi / C$$

also has nonnegative expansion coefficients by $z$ and inevitably $\exp(s(z)) - 1 = s(z) + s^2(z)/2 + \dots$ and thus $U(z)$ also have nonnegative expansion coefficients by $z \in \theta(\mathbf{C})$. But the algebra **K** is alternative and hence power associative, moreover, $\bigcup_\theta \theta(\mathbf{C}) = \mathbf{K}$. Thus expansion coefficients of $U(x_1 + x_2 + \dots + x_n)$ by degrees of $x_1,\dots,x_n$ are nonnegative, consequently, $U(0, x_1, \dots, x_n)$ is the majorant of zero. Therefore, functions $U_j(t, x_1, \dots, x_n) = U(t/\alpha + x_1 + \dots + x_n)$ compose the solution of problem II.

The problem I solved above can be reformulated in such a fashion, that there exists an inverse operator $K^{-1}$ to the differential operator together with initial conditions which compose the operator equation $Ku = g$ in the space of holomorphic vector-functions $u, g$.

Consider now the case, when $F_j$ are not **R**-homogeneous **K**-additive by $u_j$ or $v_{j,s}$ for some $j, s$. Denote the operator corresponding to this problem by $M$ and the problem is $Mu = g$ in the space of holomorphic vector-functions. Let $X$ be a **K**-vector space.

Then we call a subset $W$ in $X$ **K**-convex if $(a_1 v_1)b_1 + a_2(v_2 b_2) \in W$ for each $a_1, a_2, b_1, b_2 \in B(\mathbf{K}, 0, 1)$ and each $v_1, v_2 \in W$. Since $F_j$ are holomorphic by $(t, x)$ and locally analytic by $u_j, v_{j,s}$ for each $j, s$, then we can take power series expansions by these variables and consider a sufficiently small **K**-convex neighborhood of the point of zero initial conditions, where corrections of the second and higher order terms in $u_j, v_{j,s}$ are much smaller. This neighborhood $W$ can be chosen such that $M = K + (M - K)$ with $\|M - K\|_W / \|K\|_W < 1/2$, where $\|M\|_W := \sup_{0 \neq u \in W} |Mu|/|Ku|$. Hence $M$ is invertible on $W$ such that the inverse operator exists $M^{-1} : M(W) \to W$.

The inverse of a holomorphic function is holomorphic (see also Theorem 2). That is functions defining $M^{-1}$ are holomorphic by $(x, u)$ also by $t$, when $t \in \mathbf{K}$, moreover, they are locally analytic by others variables. In view of the proof above we can construct a solution $u$ of $Ku(t, x) = g$ for each $(t, x)$ with $|(t, x)| < R$ for some $0 < R < \infty$ such that $u \in W/2$. Therefore, we get a solution of the problem $Mu = g$ also for each $(t, x)$ in a sufficiently small neighborhood of 0. In view of the proof above this leads to the existence of a solution of the general problem.

**3. Note.** The problem 1 is the non-commutative holomorphic analog of Cauchy's problem, while Theorem 2 is the non-commutative analog of Kovalevsky's theorem (see its commutative classical form in §2 [30]). Instead of initial conditions at each $t$ with $Re(t) = Re(t_0)$ it is possible for $t$ in a $2^r - 1$-connected domain $P$ in either the quaternion skew field **H** or in the octonion algebra **O** consider initial conditions for each $t \in \partial P$, when $P$ in this section instead of $U$ there satisfies conditions of Remark 2.1 below. The proof is quite analogous in this case. This follows also by applying the non-commutative analog of the Riemann mapping Theorem 2.43 [25] establishing $\mathcal{A}_r$-holomorphic diffeomorphism between different domains in $\mathcal{A}_r$ under definite rather mild conditions. For convenience of the readers we give below the detailed proof of the chain rule over the Cayley-Dickson algebra $\mathcal{A}_r$ (though it is contained in [21, 22] as well).

**4. Notation.** If $f : U \to \mathcal{A}_r$ is either $z$-differentiable or $\tilde{z}$-differentiable at $a \in U$ or on $U$, then we can write also $D_{\tilde{z}}$ instead of $\partial_{\tilde{z}}$ and $D_z$ instead of $\partial_z$ at $a \in U$ or on $U$ respectively in situations, when it can not cause a confusion, where $U$ is open in $\mathcal{A}_r^n$.

**5. Proposition.** *Let $g : U \to \mathcal{A}_r^m$, $r \geq 2$, and $f : W \to \mathcal{A}_r^n$ be two differentiable functions on $U$ and $W$ respectively such that $g(U) \subset W$, $U$ is open in $\mathcal{A}_r^k$, $W$ is open in $\mathcal{A}_r^m$, $k, n, m \in \mathbf{N}$, where $f$ and $g$ are simultaneously either $(z, \tilde{z})$, or $z$, or $\tilde{z}$-differentiable. Then the composite function $f \circ g(z) := f(g(z))$ is differentiable on $U$ and*

$$(Df \circ g(z)).h = (Df(g)).((Dg(z)).h)$$

*for each $z \in U$ and each $h \in \mathcal{A}_r^k$, and hence $f \circ g$ is of the same type of differentiability as $f$ and $g$.*

**Proof.** Theorems 2.11, 2.15, 2.16, 3.10 and Corollary 2.13 [21, 22] establish the equivalence of notions of $\mathcal{A}_r$-holomorphic, $\mathcal{A}_r$-integral holomorphic and $\mathcal{A}_r$ locally $z$-analytic classes of functions on domains satisfying definite conditions.

In view of these results it is sufficient to prove this Proposition on open domains $U$ and $W$, where $g(U) \subset W$ with $U = g^{-1}(W)$, if others conditions are the same, since $z$-differentiability is equivalent with the local $z$-analyticity (that is in the $z$-representation). Indeed, the composition of two locally $z$-analytic functions $f \circ g$ with such domains is the locally analytic function.

Since $g$ is differentiable, the function $g$ is continuous and $g^{-1}(W)$ is open in $\mathcal{A}_r^k$. In view of Proposition 3 above if $f$ and $g$ are simultaneously either $z$-differentiable or $\tilde{z}$-differentiable, then either $\partial_{\tilde{z}} f = 0$ and $\partial_{\tilde{z}} g = 0$ or $\partial_z f = 0$ and $\partial_z g = 0$ correspondingly on their domains.

Consider the increment of the composite function

$$f \circ g(z+h) - f \circ g(z) = (Df(g))|_{g=g(z)}.(g(z+h) - g(z)) + \epsilon_f(\eta)|\eta|,$$

where $\eta = g(z+h) - g(z)$, $g(z+h) - g(z) = (Dg(z)).h + \epsilon_g(h)|h|$ (see §2.2 [21, 22]). Since the derivative $Df$ is $\mathcal{A}_r^m$-additive and $\mathbf{R}$-homogeneous (and continuous) operator on $\mathcal{A}_r^m$, we have

$f \circ g(z+h) - f \circ g(z) = (Df(g))|_{g=g(z)}.((Dg(z)).h) + \epsilon_{f \circ g}(h)|h|$, where

$$\epsilon_{f \circ g}(h)|h| := \epsilon_f((Dg(z)).h + \epsilon_g(h)|h|)|(Dg(z)).h + \epsilon_g(h)|h||$$

$$+[(Df(g))|_{g=g(z)}.(\epsilon_g(h))]|h|),$$

$|(Dg(z)).h + \epsilon_g(h)|h|| \leq [\|Dg(z)\| + |\epsilon_g(h)|]|h|$, hence

$|\epsilon_{f \circ g}(h)| \leq |\epsilon_f((Dg(z)).h + \epsilon_g(h)|h|)|[\|Dg(z)\| + |\epsilon_g(h)|] + \|(Df(g))|_{g=g(z)}\||\epsilon_g(h)|$

and inevitably $\lim_{h \to 0} \epsilon_{f \circ g}(h) = 0$. Moreover, $\epsilon_{f \circ g}(h)$ is continuous in $h$, since $\epsilon_g$ and $\epsilon_f$ are continuous functions, $Df$ and $Dg$ are continuous operators. Evidently, if $\partial_{\tilde{z}} f = 0$ and $\partial_{\tilde{z}} g = 0$ on domains of $f$ and $g$ respectively, then $\partial_{\tilde{z}} f \circ g = 0$ on $V$, since $D = \partial_z + \partial_{\tilde{z}}$.

Suppose now that there are phrases corresponding to $f_n$ and $g_n$ denoted by $\mu_n$ and $\nu_n$ such that $f_n$ and $g_n$ uniformly converge to $f$ and $g$ respectively on each bounded canonical closed subset in $W$ and $U$ from a family $\mathcal{W}$ or $\mathcal{U}$ respectively, where $\mathcal{W}$ and $\mathcal{U}$ are coverings of $W$ and $U$ correspondingly. While $D_z \mu_n$ and $D_z \nu_n$ are fundamental sequences uniformly on each bounded canonical closed subset $P \in \mathcal{W}$ and $Q \in \mathcal{U}$ correspondingly relative to the operator norm (see Definitions in §2 [21, 22]). Then the sequence $\mu_n \circ \nu_n$ converges on each bounded canonical closed subset $Q_1$ such that $Q_1 \subset Q \cap g^{-1}(P)$, where $P \in \mathcal{W}$ and $Q \in \mathcal{U}$ in $U$, when $n$ tends to the infinity. Moreover, $D_y \mu_n(y, \tilde{y}).(D_z \nu_n(z, \tilde{z}))|_{y=\nu_n(z,\tilde{z})}$ and $D_{\tilde{y}} \mu_n(y, \tilde{y}).(D_z \tilde{\nu}_n(z, \tilde{z}))|_{y=\nu_n(z,\tilde{z})}$ are the fundamental sequences of operators on each bounded canonical closed subset $Q_1$ so that $Q_1 \subset Q \cap g^{-1}(P)$, where $P \in \mathcal{W}$ and $Q \in \mathcal{U}$. The family $\mathcal{Q} := \{Q_1\}$ specified above evidently is the covering of $V$.

It remains to verify, that $D_z[\mu_n \circ \nu_n(z, \tilde{z})] = [(D_y \mu_n(y, \tilde{y}).(D_z \nu_n(z, \tilde{z})) + (D_{\tilde{y}} \mu_n(y, \tilde{y}).(D_z \tilde{\nu}_n(z, \tilde{z}))]|_{y=\nu_n(z,\tilde{z})}$.

Since the derivation operator $D_z$ by $z$ is $\mathbf{R}$-homogeneous and $\mathcal{A}_r$-additive, it is sufficient to verify this in the case $D_z[\eta \circ \psi(z, \tilde{z})]$ locally in balls, where both series uniformly converge. Here the phrase is written as:

(1) $\quad \eta = \eta(z, \tilde{z}) = \sum_k \{A_k, z, \tilde{z}\}_{q(k)} := \{a_{k,1} \widehat{z_{l_1}}^{k_1} ... a_{k,p} \widehat{z_{l_p}}^{k_p}\}_{q(k)}$,

$k = (k_1, ..., k_p)$, $p = p(k) \in \mathbf{N}$, $0 \leq k_j \in \mathbf{Z}$ and $l_j \in \{1, 2\}$ for each $j$, $a_{k,1}, ..., a_{k,p} \in \mathcal{A}_r$ are constants, $\widehat{z_1} = z$, $\widehat{z_2} = \tilde{z}$, $z$ is the Cayley-Dickson variable, a vector $q$ indicates on an order of the multiplication,

$A_k = [a_{k,1}, ..., a_{k,p}]$, also

$$(2) \quad \psi = \psi(z, \tilde{z}) = \sum_m \{B_m, z, \tilde{z}\}_{q(m)},$$

where $m = (m_1, ..., m_s)$, $B_m = [b_{m,1}, ..., b_{m,s}]$, $b_{m,j} \in \mathcal{A}_r$ is a constant for each $m, j$, also $s \in \mathbf{N}$, $0 \le m_j \in \mathbf{N}$ for each $j$. We have the identities

$$(3) \quad D_z \widehat{z_l} = \delta_{l,1} \mathbf{1}, \quad D_{\tilde{z}} \widehat{z_l} = \delta_{l,2} \tilde{\mathbf{1}},$$

where $\mathbf{1}h := h$ and $\tilde{\mathbf{1}}h := \tilde{h}$ for each $h \in \mathcal{A}_r$. Using shifts $z \mapsto z - \zeta$ we can consider that the series are decomposed around a zero point. If a series $\eta$ is uniformly converging on a canonical closed subset $Q_1$ as above, then

$$(4) \quad D_z \eta = \sum_k D_z \{a_{k,1} \widehat{z_{l_1}^{k_1}} ... a_{k,p} \widehat{z_{l_p}^{k_p}}\}_{q(k)}$$

is also uniformly converging on $Q_1$ in accordance with our supposition about convergence of phrases and their derivatives.

In suitable $Q_1$ we deduce from Formula (4) that

$$(5) \quad D_z[\eta \circ \psi(z, \tilde{z})] = D_z \sum_k \{a_{k,1} \widehat{\psi_{l_1}^{k_1}} ... a_{k,p} \widehat{\psi_{l_p}^{k_p}}\}_{q(k)}$$

$$= \sum_k D_z \{a_{k,1} \widehat{\psi_{l_1}^{k_1}} ... a_{k,p} \widehat{\psi_{l_p}^{k_p}}\}_{q(k)}$$

$$= \sum_k [\{a_{k,1}(D_z \widehat{\psi_{l_1}^{k_1}}) ... a_{k,p} \widehat{\psi_{l_p}^{k_p}}\}_{q(k)} + ... + \{a_{k,1} \widehat{\psi_{l_1}^{k_1}} ... a_{k,p}(D_z \widehat{\psi_{l_p}^{k_p}})\}_{q(k)}]$$

$$= \sum_k \sum_{{}_1m, ..., {}_pm; p=p(k)} [\{a_{k,1}(D_z(\widehat{\{B_{{}_1m}, z, \tilde{z}\}_{q({}_1m)}})^{k_1}_{l_1}) ... a_{k,p}(\widehat{\{B_{{}_pm}, z, \tilde{z}\}_{q({}_pm)}})^{k_p}_{l_p}\}_{q(k)}$$

$$+ ... + \{a_{k,1}(\widehat{\{B_{{}_1m}, z, \tilde{z}\}_{q({}_1m)}})^{k_1}_{l_1} ... a_{k,p}(D_z(\widehat{\{B_{{}_pm}, z, \tilde{z}\}_{q({}_pm)}})^{k_p}_{l_p})\}_{q(k)}.$$

On the other hand, we have

$$D_z\{B_m, z, \tilde{z}\}_{q(m)} = \{b_{m,1}(D_z \widehat{z_{l_1}^{m_1}}) ... b_{m,s} \widehat{z_{l_s}^{k_s}}\}_{q(m)} + ... + \{b_{m,1} \widehat{z_{l_1}^{m_1}} ... b_{m,s}(D_z \widehat{z_{l_s}^{k_s}})\}_{q(k)}$$

due to the Leibnitz rule, hence

$$D_z \widehat{\psi^p} = \sum_{{}_1m, ..., {}_pm} [(...(D_z \widehat{\{B_{{}_1m}, z, \tilde{z}\}_{q({}_1m)}})...) \widehat{\{B_{{}_pm}, z, \tilde{z}\}_{q({}_pm)}}$$

$$+ ... + (...(\widehat{\{B_{{}_1m}, z, \tilde{z}\}_{q({}_1m)}} \widehat{\{B_{{}_2m}, z, \tilde{z}\}_{q({}_2m)}})...)(D_z \widehat{\{B_{{}_pm}, z, \tilde{z}\}_{q({}_pm)}})]$$

$$= (D_y \widehat{y^p}).(D_z \psi) + (D_{\tilde{y}} \widehat{y^p}).(D_z \tilde{\psi})$$

due to Formulas (3). Thus Formulas (5, 6) imply that

(7) $D_z(\eta\circ\psi(z,\tilde{z})) = [(D_y\eta(y,\tilde{y})).(D_z\psi(z,\tilde{z}))+(D_{\tilde{y}}\eta(y,\tilde{y})).(D_z\tilde{\psi}(z,\tilde{z}))]|_{y=\psi(z,\tilde{z})}.$

Quite analogously or using the conjugation one deduces that

(8) $D_{\tilde{z}}(\eta\circ\psi(z,\tilde{z})) = [(D_{\tilde{y}}\eta(y,\tilde{y})).(D_{\tilde{z}}\psi(z,\tilde{z}))+(D_y\eta(y,\tilde{y})).(D_{\tilde{z}}\tilde{\psi}(z,\tilde{z}))]|_{y=\psi(z,\tilde{z})}.$

Particularly, we get

(9) $D_z(\eta\circ\psi(z)) = (D_y\eta(y)).(D_z\psi(z))|_{y=\psi(z)}$, when $D_{\tilde{z}}\eta(z,\tilde{z}) = 0$ and $D_{\tilde{z}}\psi(z,\tilde{z}) = 0$

and

(10) $D_{\tilde{z}}(\eta\circ\psi(\tilde{z})) = (D_{\tilde{y}}\eta(\tilde{y})).(D_{\tilde{z}}\psi(\tilde{z}))|_{y=\psi(\tilde{z})}$, when $D_z\eta(z,\tilde{z}) = 0$ and $D_z\psi(z,\tilde{z}) = 0$.

Combining Formulas (7, 8) and taking into account the proof given above and applying Proposition 2.3 [21, 22] we infer the chain rule, when both $\eta$ and $\psi$ are either the $z$ or $\tilde{z}$ or $(z,\tilde{z})$-differentiable, that is in all three considered cases.

**6. Proposition.** *Let suppositions of Proposition 5 be satisfied for two ($\mathcal{A}_r$ super-) $z$-differentiable functions $f$ and $g$, where $U$ and $W \subset \mathcal{A}_r$. Then*

(1) $(d^n f(g(z))/dz^n).(h_1, ..., h_n) = f^{(n)}(g(z)).(g'(z).h_1, ..., g'(z).h_n) +$

$\sum_{k,j=1,...,n;\ k\neq j;\ j_1<...<j_{n-2}\in\{1,...,n\}\setminus\{k,j\}} f^{(n-1)}(g(z)).(g^{(2)}.(h_k, h_j), g'(z).h_{j_1}, ..., g'(z).h_{j_{n-2}}) +$

$\sum_{2\leq k\leq n-2;\ p_1+...+p_l=n,\ p_1\geq p_2\geq...\geq p_l\geq 1,\ l=n-k}$

$f^{(n-k)}(g(z)).(g^{(p_1)}.(h_{j(1,p_1)}, ..., h_{j(p_1,p_1)}), ..., g^{(p_l)}.(h_{j(1,p_l)}, ..., h_{j(p_l,p_l)}))$

$+f'(g(z)).(g^{(n)}.(h_1, ..., h_n))$

*for all Cayley-Dickson numbers $h_1, ..., h_n \in \mathcal{A}_r$, $2 \leq r$, $n \geq 2$, where $j(s,p) \in \{1, ..., n\}$ are all pairwise distinct natural numbers for different $(s,p)$, $j(s,p_k) < j(s+1,p_k)$ for all $s$ and $k$, $j(1,k) < j(1,k+1)$ if $p_k = p_{k+1}$.*

**Proof.** In accordance with §5 it is sufficient to prove Formula (1) for phrases $\eta(z)$ and $\psi(z)$, since each $z$-differentiable function is infinite $z$-differentiable (see §2 [21, 22]). To derive Formula (1) the mathematical induction will be used below. For $n = 1$ Formula (1) coincides with that of Proposition 4. In accordance with the Leibnitz rule over the Cayley-Dickson algebra $\mathcal{A}_r$, $2 \leq r$, [21, 22] the Formula

(2) $(d\{\eta_1(z)...\eta_n(z)\}_q/dz).h = \sum_{j=1}^{n}\{\eta_1...[(d\eta_j(z)/dz).h]...\eta_n\}_q$

is accomplished, where $\eta_1,...,\eta_n$ are phrases of $z$-differentiable functions on the domain $W$ open in $\mathcal{A}_r$. Mention that the differential operator $\eta^{(k)}(z).(h_1,...,h_k)$ is as usually symmetric relative to all transpositions of $h_1,...,h_n$ (see also [21, 22]). Therefore, we deduce from the Leibnitz rule, that

(3) $(d^2\eta(\psi(z))/dz^2).(h_1,h_2) = d[\eta'(\psi(z)).[\psi'(z).h_1]/dz].h_2 = \eta^{(2)}(\psi(z)).([\psi'(z).h_1],[\psi'(z).h_2]) + \eta'(\psi(z)).[\psi^{(2)}(z).(h_1,h_2)]$ and

(4) $[d^3\eta(\psi(z))/dz^3].(h_1,h_2,h_3) = d\{[d^2\eta(\psi(z))/dz^2].(h_1,h_2)\}.h_3$
$= \eta^{(3)}(\psi(z)).([\psi'(z).h_1],[\psi'(z).h_2],[\psi'(z).h_3]) + \eta^{(2)}(\psi(z)).([\psi^{(2)}(z).(h_1,h_3)],[\psi'(z).h_2]) + \eta^{(2)}(\psi(z)).([\psi^{(2)}(z).(h_2,h_3)],[\psi'(z).h_1]) + \eta^{(2)}(\psi(z)).([\psi^{(2)}(z).(h_1,h_2)],[\psi'(z).h_3]) + \eta'(\psi(z)).[\psi^{(3)}(z).(h_1,h_2,h_3)]$.

Suppose now that Formula (1) for $d^k\eta(\psi(z))/dz^k$ is proved for each $k = 1,...,n$. It is sufficient according to the continuation of the induction to demonstrate it for $k = n+1$. From $(1,2)$ we infer the equality:

(5) $[d^{n+1}\eta(\psi(z))/dz^{n+1}].(h_1,...,h_{n+1}) = \{d[d^n\eta(\psi(z))/dz^n].(h_1,...,h_n)/dz\}.h_{n+1} = \eta^{(n+1)}(\psi(z)).(\psi'(z).h_1,...,\psi'(z).h_{n+1}) +$

$\sum_{j=1,...,n;\ 1\leq j_1<...<j_{n-1}\leq n,\ j_k\neq j \forall k} \eta^{(n)}(\psi(z)).([\psi^{(2)}(z).(h_j,h_{n+1})],[\psi'(z).h_{j_1}],...,\psi'(z).h_{j_{n-1}}) +$

$\sum_{k,j=1,...,n;\ k\neq j;\ j_1<...<j_{n-2}\in\{1,...,n\}\setminus\{k,j\}}$
$\eta^{(n)}(\psi(z)).(\psi^{(2)}.(h_k,h_j),\psi'(z).h_{j_1},...,\psi'(z).h_{j_{n-2}},\psi'(z).h_{n+1}) +$

$\sum_{2\leq k\leq n-1;\ p_1+...+p_l=n,\ p_1\geq p_2\geq...\geq p_l\geq 1,\ l=n-k, p_{l+1}=1} \eta^{(n+1-k)}$
$(\psi(z)).(\psi^{(p_1)}.(h_{j(1,p_1)},...,h_{j(p_1,p_1)}),...,\psi^{(p_l)}.(h_{j(1,p_l)},...,h_{j(p_l,p_l)}),\psi'(z).h_{n+1})$

$\sum_{2\leq k\leq n-1;\ p_1+...+p_l=n,\ p_1\geq p_2\geq...\geq p_l\geq 1,\ l=n-k,\ s=1,...,l}$
$\eta^{(n+1-k)}(\psi(z)).(\psi^{(p_1)}.(h_{j(1,p_1)},...,h_{j(p_1,p_1)}),...,\psi^{(p_s+1)}.(h_{j(1,p_s)},...,h_{j(p_s,p_s)},h_{n+1}),...,$
$\psi^{(p_l)}.(h_{j(1,p_l)},...,h_{j(p_l,p_l)})) + \eta'(\psi(z)).(\psi^{(n+1)}.(h_1,...,h_{n+1}))$

for all Cayley-Dickson numbers $h_1,...,h_{n+1}$. This is the same expression as in Formula (1) for $n+1$ instead of $n$ after gathering together terms $\eta^{(j)}$ with the same $j$.

**7. Note.** In Propositions 5 and 6 the specification of phrases corresponding to functions in the respective representation is essential as well as for calculations of their derivative operators. Finally one can mention, that for restrictions of $f$ and $g$ on an open domain $U_\mathbf{C}$ in the complex field $\mathbf{C}$ embedded as $\mathbf{C_M} := \mathbf{R}i_0 \oplus M\mathbf{R}$ into the Cayley-Dickson algebra $\mathcal{A}_r$ for some purely imaginary Cayley-Dickson number $M$ of unit absolute value and for $h_1 = ..., h_n = 1$ when $f(z), g(z) \in \mathbf{C}$ for each $z \in U_{\mathbf{C}}$ the deduced formula coincides with that of 5.2(1) in [32].

In this article some elementary facts about analytic functions $z^n$, $z^{1/n}$, $\exp(z)$ and $Ln(z)$ of the Cayley-Dickson variables are used. They were considered in details in previous works [22, 21, 25]. Recall that the exponential function is defined by the power series $\exp(z) := 1 + \sum_{n=1}^{\infty} z^n/n!$ converging on the entire Cayley-Dickson algebra $\mathcal{A}_r$. It has the periodicity property $\exp(M(\phi+2\pi k)) = \exp(M\phi)$ for each purely imaginary Cayley-Dickson number $M$ of the unit norm $|M| = 1$ for any real number $\phi \in \mathbf{R}$ and each integer number $k \in \mathbf{Z}$. The restriction of such exponential function on each complex plane $\mathbf{C}_M := \mathbf{R} \oplus M\mathbf{R}$ coincides with the traditional complex exponential function. Since the inverse function $Ln(z)$ of $z = \exp(x)$ is defined on every complex plane $\mathbf{C}_M \setminus \{0\}$ with the pricked zero point, the logarithmic function $Ln(z)$ is defined on $\mathcal{A}_r \setminus \{0\}$. This logarithmic function is certainly multivalued. Consider the bunch of complex planes $\mathbf{C}_M$ intersecting by the real line $\mathbf{R}i_0$ as the geometric realization of $\mathcal{A}_r$. Certainly $\mathbf{C}_M = \mathbf{C}_{-M}$, so we take the set

$\mathbf{S}_r^+ := \{M \in \mathcal{A}_r : |M| = 1, M = M_1 i_1 + ... + M_{2^r-1} i_{2^r-1}$, either $M_1 > 0$, or $M_1 = 0$ and $M_2 > 0$, or..., or $M_1 = ... = M_{2^r-2} = 0$ and $M_{2^r-1} > 0\}$, where $M_1, ..., M_{2^r-1} \in \mathbf{R}$, $2 \le r$. Then $\mathcal{A}_r = \bigcup_{M \in \mathbf{S}_r^+} \mathbf{C}_M$.

If $A$ and $B$ are two subsets in a complete uniform space $X$ and $\theta : A \to B$ is a continuous bijective mapping, the equivalence relation $a\Upsilon b$ by definition means $b = \theta(a)$ for $a \in A$ and $b \in B$; or $b\Upsilon a$ means $a = \theta^{-1}(b)$. When $\theta$ is uniformly continuous, $\theta$ has a uniformly continuous extension $\theta : cl(A) \to cl(B)$, where $cl(A)$ denotes the closure of the set $A$ in $X$ (see Theorem 8.3.10 in [5]). Certainly, the mapping $\theta$ can be specified by its graph $\{(x_1, x_2) : x_2 = \theta(x_1), x_1 \in A\}$. We say, that $A$ and $B$ are glued (by $\theta$), if $B = \theta(A)$ and the natural quotient mapping $\pi : X \to X/\Upsilon$ is given, where $X/\Upsilon$ denotes the quotient space (see §2.4 [5]). For the uniformly continuous $\theta$ this means that the gluing is extended from $A$ onto $cl(A)$.

In the complex case to construct the Riemann surface of the logarithmic function one takes traditionally the complex plane $\mathbf{C}$ cut by the set $Q_1 := \{z = x + iy \in \mathbf{C} : x < 0\}$ and marking two respective points $x_1$ and $x_2$ of two edges $Q_{1,1}$ and $Q_{1,2}$ of the cut $Q_1$ arising from each given point $x < 0$, where $i = i_1$. Then one embeds $\mathbf{C}$ into either $\mathbf{C} \times \mathbf{R}$ or $\mathbf{C} \times \mathbf{R}i$ and bents the obtained surface slightly along the perpendicular axis $e_3$ to $\mathbf{C}$ by neighborhoods of two edges of the cut $Q_1$ and gets the new surface $\mathcal{C}$. Taking the countable infinite family of such surfaces $\mathcal{C}^j$ with the edges of cuts $Q_{1,1}^j$ and $Q_{1,2}^j$ and gluing by respective points of the cuts $Q_{1,2}^j$ with $Q_{1,1}^{j+1}$ for each $j$ one gets the Riemann surface of the logarithmic function, where $j \in \mathbf{Z}$ (see, for example, [18]).

Analogous procedure to construct the Riemann surface is in the cases $r \geq 2$: one cuts $\mathcal{A}_r$ by $Q_r$ and gets two edges $Q_{r,1}$ and $Q_{r,2}$ of the cut. This is described below.

If $K$ and $M = KL$ are two purely imaginary Cayley-Dickson numbers with $|K| = |M| = |L| = 1$ so that they are orthogonal $K \perp M$, that is $Re(KM) = 0$, then $(K\mathbf{R}) \oplus (M\mathbf{R}) = K(\mathbf{R} \oplus L\mathbf{R})$ and $L$ is also purely imaginary. Consider any path $\gamma : [0,1] \to K\mathbf{R} \oplus M\mathbf{R}$ winding one time around zero such that $\gamma(t) \neq 0$ for each $t$. Each $z \in K\mathbf{R} \oplus M\mathbf{R}$ can be written in the polar form $z = |z|Ke^{L\phi} = |z|\exp(\pi Ke^{L\phi}/2)$, where $\phi = \phi(z) \in \mathbf{R}$, since $Ke^{L\phi} = K\cos(\phi) + (KL)\sin(\phi)$ due to Euler's formula, hence $|Ke^{L\phi(\gamma(t))}| = 1$ for each $t$ (see Section 3 in [22, 21]). In particular, $\gamma(t) = |\gamma(t)|\exp(\pi Ke^{L\phi(t)}/2)$. But $Ke^{L\phi(\gamma(0))} = Ke^{L\phi(\gamma(1))}$ such that the logarithm $Ln\ \gamma(t) = Ln\ |\gamma(t)| + \pi Ke^{L\phi(\gamma(t))}/2$ does not change its branch, when the path $\gamma(t) \in (K\mathbf{R}) \oplus (M\mathbf{R}) \setminus \{0\}$ winds around zero. Due to the homotopy theorem (see [22, 21, 26]) this means that the logarithm $Ln\ \gamma(t)$ does not change its branch, when $Re(\gamma(t)) = 0$ for each $t$ for the path $\gamma$ winding around zero with $|\gamma(t)| > 0$ for each $t$.

The first simple construction for $r \geq 2$ is the following. Take the set $Q_r := (-\infty, 0)\mathbf{S}_r^+ := \{z = tx : t \in (-\infty, 0), x \in \mathbf{S}_r^+\}$ and cut $\mathcal{A}_r$ by $Q_r$. The set $Q_r$ is the union $Q_r = \bigcup_{j=1}^{2^r-1} \Omega_{j,r}$ of subsets $\Omega_{j,r} := \{z \in Q_r : z_0 = 0, ..., z_j = 0\}$ so that $\Omega_{j,r}$ is contained in the boundary of the preceding set $\Omega_{j,r} \subset \partial \Omega_{j-1,r}$ for each $j = 1, ..., 2^r - 1$, $dim\ \Omega_{j,r} = dim\ \Omega_{j-1,r} - 1$, moreover, $\mathbf{RS}_r^+ = \mathcal{I}_r := \{z \in \mathcal{A}_r : Re(z) = 0\}$. Therefore, from each point $z \in Q_r$ two and only two different points $z_1$ and $z_2$ arise.

It is useful to embed $\mathcal{A}_r$ either into $\mathcal{A}_r \times \mathbf{R}^{2^r-1}$ or into $\mathcal{A}_r \times \mathcal{I}_r$, where $\mathcal{I}_r := \{z \in \mathcal{A}_r : Re(z) = 0\}$. Then one marks all pairs of respective points $z_1$ and $z_2$ arising from $z \in Q_r$ after cutting, slightly bents the cut copy of $\mathcal{A}_r \setminus \{0\}$ by $(2^r - 1)$ axes perpendicular to $\mathcal{A}_r$ by two neighborhoods of two edges of the cut $Q_r$. Thus one gets the $2^r$ dimensional surface $\mathcal{C}_r$ with two edges $Q_{r,1}$ and $Q_{r,2}$ of the cut. Taking the countable infinite family of such surfaces $\mathcal{C}_r^j$ with edges of the cuts $Q_{r,1}^j$ and $Q_{r,2}^j$, $j \in \mathbf{Z}$, and gluing respective points of edges $Q_{r,2}^j$ with $Q_{r,1}^{j+1}$ for each $j$ one gets the Riemann surface $\mathcal{R}_r = \mathcal{R}_{r,Ln}$ of the logarithmic function $Ln(z) : \mathcal{A}_r \setminus \{0\} \to \mathcal{R}_r$. Thus the latter mapping is already univalent with the image in the Riemann surface (see in details [22, 21, 25]).

For convenience we attach numbers 1 and 2 to faces in such manner that the winding around zero in the complex plane $\mathbf{C}_M$ embedded into $\mathcal{R}_r$ counterclockwise means the transition through the cut from $Q_{r,2}^j$ to $Q_{r,1}^j$ for each $M$ in the connected set $\mathbf{S}_r^+$.

For the function $z^{1/n}$ with $n \in \mathbf{N}$ its Riemann surface $\mathcal{R}_{r,z^{1/n}}$ is obtained from $n$ copies of surfaces $\mathcal{C}_r^j$, $j = 1, ..., n$, by gluing the corresponding points of edges $Q_{r,2}^j$ with $Q_{r,1}^{j+1}$ for $j = 1, ..., n-1$ and of $Q_{r,2}^n$ with $Q_{r,1}^1$ (see also [18, 22, 21, 25]).

Another more complicated construction is described below. Now we take the set

$Q_r := \bigcup_{j=1}^{2^r-1} P_j$, where

$P_j := \{z \in \mathcal{A}_r : z = z_0 i_0 + ... + z_{2^r-1} i_{2^r-1}; z_0 < 0 \text{ and } z_j = 0\}$, where $z_0, ..., z_{2^r-1} \in \mathbf{R}$, $2 \le r$.

Let $z = z_0 + z'$ be the Cayley-Dickson number with the negative real part $z_0 < 0$ and the imaginary part $Im(z) = z'$, which can be written in the form $z' = |z'| \exp(\pi K e^{L\phi(z')}/2)$ and $z = |z| e^{P\psi}$, where $K$, $L$ and $P$ are purely imaginary of the unit norm, $\phi$ and $\psi \in \mathbf{R}$ are reals. This gives the relations $\cos(\psi) = z_0/|z|$ so that the parameter $\psi$ is in the interval $\pi/2 + 2\pi k < \psi < 3\pi/2 + 2\pi k$ for some integer number $k$. This means that for a continuous path $\gamma$ contained in the set $Q_r$ the parameter $\psi(\gamma(t))$ is the continuous function of the real variable $t \in \mathbf{R}$ and remains in the same interval $(\pi/2 + 2\pi k, 3\pi/2 + 2\pi k)$, consequently, $Ln\gamma(t)$ preserves it branch along such path $\gamma(t)$. This shows that after the first cut along $Q_r$ the obtained sets $Q_{r,1}$ and $Q_{r,2}$ need not be further cut. Thus the described reason simplifies the construction of the Riemann surface.

Each continuous path $\gamma : [0,1] \to \mathcal{A}_r$ can be decomposed as the sum as well as the composition up to the homotopy satisfying the conditions of the homotopy theorem [22, 21, 26] of paths $\gamma_{k,l}$ in complex planes $(\mathbf{R}i_k) \oplus (\mathbf{R}i_l)$ for each $k < l \in \Lambda$ for the corresponding subset $\Lambda \subset \{0, 1, ..., 2^r - 1\}$. If $|\gamma(t)| > 0$ for each $t$ we take $\gamma_{k,l}$ with $|\gamma_{k,l}(t)| > 0$ on $[0,1]$ for all $k < l \in \Lambda$. When $\gamma[0,1]$ does not intersect $Q_r$ one can choose $\gamma_{k,l}$ with images $\gamma_{k,l}[0,1]$ also non-intersecting with $Q_r$ for all $k < l \in \Lambda$. Due to the homotopy theorem the logarithm $Ln\, \gamma(t)$ does not change its branch along such continuous path $\gamma(t)$, since this is the case for $Ln\, \gamma_{k,l}(t)$ for all $k < l \in \Lambda$.

For each $z \in P_j \setminus \bigcup_{m, m \neq j} P_m$ cutting by $Q_r$ gives two points. If $z \in (P_k \cap P_j) \setminus \bigcup_{m, m \neq k, m \neq j} P_m$ with $k < j$ cutting gives four points which can be organized into respective pairs after cutting of $P_k$ and then of $P_j$. This procedure gives pairs $(z_{1,1}; z_{2,2})$ and $(z_{1,2}; z_{2,1})$. For $z \in (P_{k_1} \cap ... \cap P_{k_l}) \setminus \bigcup_{m, m \neq k_1, ..., m \neq k_l} P_m$ consider pairs of points appearing from the preceding point after cutting of $P_{k_j}$, $j = 1, ..., l$ by induction, where $k_1 < ... < k_l$. One can do it by induction by all $l = 1, ..., 2^r - 1$ and all possible subsets $1 \leq k_1 < k_2 < ... < k_l \leq 2^r - 1$. It produces points denoted by $z_{a_1, ..., a_l}$. Here $a_j = 1$ corresponds to the face of the cut with $z_j \leq 0$, while $a_j = 2$ corresponds to the face of the cut with $z_j \geq 0$. Two points $z_{a_1, ..., a_l}$ and $z_{b_1, ..., b_l}$ form the pair of respective points, when $a_j + b_j = 3$ for each $j$, where $l \in \{1, ..., 2^r - 1\}$ is the index of such points. To each point $z \in (P_{k_1} \cap ... \cap P_{k_l}) \setminus \bigcup_{m, m \neq k_1, ..., m \neq k_l} P_m$ of index $l \geq 1$ the point $M \in \mathbf{S}_r^+$ corresponds such that $M = M_1 i_1 + ... + M_{2^r - 1} i_{2^r - 1}$ and $z \in \mathbf{C}_M$, where $M_{k_1} = 0, ..., M_{k_l} = 0$.

Then as above after cutting of $\mathcal{A}_r$ by $Q_r$ one gets the $2^r$ dimensional surface $\mathcal{C}_r$ with two edges $Q_{r,1}$ and $Q_{r,2}$ of the cut. For the countable infinite family of the surfaces $\mathcal{C}_r^j$ with edges of the cuts $Q_{r,1}^j$ and $Q_{r,2}^j$, $j \in \mathbf{Z}$, one glues respective points of edges $Q_{r,2}^j$ with $Q_{r,1}^{j+1}$ for each $j$. Thus one gets the Riemann surface $\mathcal{R}_r = \mathcal{R}_{r,Ln}$ of the logarithmic function $Ln(z) : \mathcal{A}_r \setminus \{0\} \to \mathcal{R}_r$. Therefore, the latter mapping is univalent with the image in the Riemann surface (see in details [22, 21, 25]). For simplicity one can glue at first pairs of respective points of index $l = 1$ and extend gluing by continuity on points of index $l > 1$ (see above). The first construction for $r \geq 2$ operates with points of index one only. That is why it is simpler than the second procedure.

The reader can lightly see that geometrically $\mathcal{R}_{r,f}$ (for both types described above) is the bunch $\bigcup \{\mathcal{R}_{1,f,M} : M \in \mathbf{S}_r^+\}$ of complex Riemann surfaces $\mathcal{R}_{1,f,M}$ for the restrictions $f|_{\mathbf{C}_M}$ on the complex planes $\mathbf{C}_M$ of the function $f$, where either $f = Ln$ or $f = z^{1/n}$ among those considered here.

The main problem of the multi-dimensional geometry is in depicting its objects on the two-dimensional sheet of paper so one uses either projections or sections of the multi-dimensional surfaces. The bunch interpretation $\bigcup_{M \in \mathbf{S}_r^+} \mathcal{R}_{1,f,M}$ helps to imagine how the $2^r$-dimensional Riemann surface $\mathcal{R}_{r,f}$ is organized for $f = Ln$ and $f = z^{1/n}$ with $n \in \mathbf{N}$.

## 3 Differential equations over the octonion algebra.

**1. Remarks.** For a subset $U$ in either the quaternion skew field $\mathbf{H} = \mathcal{A}_2$ or in the octonion algebra $\mathbf{O} = \mathcal{A}_3$ we put $\pi_{\mathsf{s},\mathsf{p},\mathsf{t}}(U) := \{\mathsf{u} : z \in U, z = \sum_{\mathsf{v} \in \mathbf{b}} w_\mathsf{v} \mathsf{v}, \mathsf{u} = w_\mathsf{s} \mathsf{s} + w_\mathsf{p} \mathsf{p}\}$ for each $\mathsf{s} \neq \mathsf{p} \in \mathbf{b}$, where $\mathsf{t} := \sum_{\mathsf{v} \in \mathbf{b} \setminus \{\mathsf{s},\mathsf{p}\}} w_\mathsf{v} \mathsf{v} \in \mathcal{A}_{r,\mathsf{s},\mathsf{p}} := \{z \in \mathcal{A}_r : z = \sum_{\mathsf{v} \in \mathbf{b}} w_\mathsf{v} \mathsf{v}, w_\mathsf{s} = w_\mathsf{p} = 0, w_\mathsf{v} \in \mathbf{R} \; \forall \mathsf{v} \in \mathbf{b}\}$, where $\mathbf{b} := \{i_0, i_1, ..., i_{2^r-1}\}$ is the family of standard generators of the algebra $\mathcal{A}_r$. Geometrically the domain $\pi_{\mathsf{s},\mathsf{p},\mathsf{t}}(U)$ means the projection on the complex plane $\mathbf{C}_{\mathsf{s},\mathsf{p}}$ of the intersection $U$ with the plane $\tilde{\pi}_{\mathsf{s},\mathsf{p},\mathsf{t}} \ni \mathsf{t}$, $\mathbf{C}_{\mathsf{s},\mathsf{p}} := \{a\mathsf{s} + b\mathsf{p} : a, b \in \mathbf{R}\}$, since $\mathsf{sp}^* \in \hat{b} := \mathbf{b} \setminus \{1\}$. Recall that in §§2.5-7 [21] for each continuous function $f : U \to \mathcal{A}_r$ it was defined the operator $\hat{f}$ by each variable $z \in \mathcal{A}_r$. If a function $f$ is $z$-differentiable by the Cayley-Dickson variable $z \in U \subset \mathcal{A}_r$, $2 \leq r$, then $\hat{f}(z) = dg(z)/dz$, where $(dg(z)/dz).1 = f(z)$.

A Hausdorff topological space $X$ is said to be $n$-connected for $n \geq 0$ if each continuous map $f : S^k \to X$ from the $k$-dimensional real unit sphere into $X$ has a continuous extension over $\mathbf{R}^{k+1}$ for each $k \leq n$ (see also [38]). A 1-connected space is also said to be simply connected.

It is supposed further, that a domain $U$ in $\mathcal{A}_r^m$ has the property that

($D1$) each projection $\mathbf{p}_j(U) =: U_j$ is $(2^r - 1)$-connected;

(D2) $\pi_{\mathsf{s},\mathsf{p},\mathsf{t}}(U_j)$ is simply connected in $\mathbf{C}$ for each $k = 0, 1, ..., 2^{r-1}$, $\mathsf{s} = i_{2k}$, $\mathsf{p} = i_{2k+1}$, $\mathsf{t} \in \mathcal{A}_{r,\mathsf{s},\mathsf{p}}$ and $u \in \mathbf{C}_{\mathsf{s},\mathsf{p}}$, for which there exists $z = \mathsf{u} + \mathsf{t} \in U_j$, where $e_j = (0, ..., 0, 1, 0, ..., 0) \in \mathcal{A}_r^m$ is the vector with 1 on the $j$-th place, $\mathbf{p}_j(z) = {}^jz$ for each $z \in \mathcal{A}_r^m$, $z = \sum_{j=1}^m {}^jze_j$, ${}^jz \in \mathcal{A}_r$ for each $j = 1, ..., m$, $m \in \mathbf{N} := \{1, 2, 3, ...\}$. Frequently we take $m = 1$. Henceforward we consider a domain $U$ satisfying Conditions $(D1, D2)$ if any other is not outlined.

The family of all $\mathcal{A}_r$ locally analytic functions $f(x)$ on $U$ with values in $\mathcal{A}_r$ is denoted by $\mathcal{H}(U, \mathcal{A}_r)$. It is supposed that a locally analytic function $f(x)$ is written in the $x$-representation $\nu(x)$, also denoted by $\nu = \nu^f$. The latter is equivalent to the super-differentiability of $f$ (see [22, 21, 23]). Each such $f$ is supposed to be specified by its phrase $\nu$.

For each super-differentiable function $f(x)$ its non-commutative line integral $\int_\gamma f(x)dx$ in $U$ is defined along a rectifiable path $\gamma$ in $U$. It is the integral of a differential form $\hat{f}(x).dx$, where

(I1) $\hat{f}(x) = dg(x)/dx$,

(I2) $[dg(x)/dx].1 = f(x)$ for each $x \in U$. A branch of the non-commutative line integral can be specified with the help of either the left or right algorithm (see [22, 21, 23]). We take further for definiteness the left algorithm if something another will not described. For $f \in \mathcal{H}(U, \mathcal{A}_r)$ and a rectifiable path $\gamma : [a, b] \to \mathcal{A}_r$ the integral $\int_\gamma f(x)dx$ depends only on an initial $\alpha = \gamma(a)$ and final $\beta = \gamma(b)$ points due to the non-commutative analog of the homotopy theorem in $U$, where $a < b \in \mathbf{R}$. When initial and final points or a path are not marked we denote the operation of the non-commutative line integration in the domain $U$ simply by $\int f(x)dx$ analogously to the indefinite integral.

To rewrite a function from real variables $z_j$ in the $z$-representation the following identities are used:

(1) $z_j = (-zi_j + i_j(2^r-2)^{-1}\{-z + \sum_{k=1}^{2^r-1} i_k(zi_k^*)\})/2$ for each $j = 1, 2, ..., 2^r - 1$,

(2) $z_0 = (z + (2^r - 2)^{-1}\{-z + \sum_{k=1}^{2^r-1} i_k(zi_k^*)\})/2$,

where $2 \le r \in \mathbf{N}$, $z$ is a Cayley-Dickson number decomposed as

(3) $z = z_0 i_0 + ... + z_{2^r-1} i_{2^r-1} \in \mathcal{A}_r$, $z_j \in \mathbf{R}$ for each $j$, $i_k^* = \tilde{i}_k = -i_k$ for each $k > 0$, $i_0 = 1$.

## 2. The simplest differential equation of the first order

has the form:

(1) $[dy(x)/dx].h = f(x)$,

where $h = h(x) \in \mathcal{A}_r$, $r \geq 2$, $f, h \in \mathcal{H}(U, \mathcal{A}_r)$, $|h(x)| > 0$ for $\lambda$-almost all $x \in U$, an open domain $U$ satisfies Conditions $1(D1, D2)$. Here $\lambda$ denotes the Lebesgue measure on $\mathcal{A}_r$ induced from that of on the real shadow $\mathbf{R}^{2^r}$. Using a non-commutative line integral in $U$ we get:

$$(2) \quad \int_\gamma ([dy(x)/dx].h)dx = \int_\gamma f(x)dx.$$

As usually

(3) $Re(\beta) = (\beta + \beta^*)/2$ denotes the real part of $\beta \in \mathcal{A}_r$, $\beta^* = \tilde{\beta}$ denotes the conjugated number $\beta$, so that its square norm is given by the equation

(4) $|\beta|^2 = \beta\beta^*$.

Suppose that $h = h(t)$ is a (super-)differentiable function on a domain $U$, $h: U \to \mathcal{A}_r$. Take a marked point $\alpha$ in $U$.

Let us consider the differential equation

(5) $(d\omega(x)/dx).1 = h(\omega)$, $\omega = \omega_h$.

This equation is equivalent to

(5.1) $(1/h(\omega))(d\omega(x)/dx).1 = 1$.

Integrating both sides of (5.1) by $x$ we deduce that

(5.2) $p(\omega) := \int_\gamma [(1/h(\omega))(d\omega(t)/dt).1]dt =$
$\int (1/h(\omega))d\omega = \int_\alpha^x dt = x - \alpha$,

since $(d[\int (d\omega/dx).1dx]/dx).1 = (d\omega/dx).1$,

where $\gamma(a) = \alpha$ and $\gamma(b) = x$, $\gamma: [a,b] \to U$ is a rectifiable path. Therefore,

(6) $\omega_h(x) = \omega(x) = p^{-1}(x - \alpha) + \phi_1(x^\diamond, \alpha^\diamond)$

satisfies Equation (5), where $p^{-1}$ denotes the inverse function. This function $\omega(x)$ is defined up to a summand $\phi_1(x^\diamond, \alpha^\diamond)$, where $\phi_1(x^\diamond, \alpha^\diamond)$ is any (super-)differentiable function of the arguments

(7) $x^\diamond := x^\diamond_1 := x - Re(x) =: Im(x)$ and $\alpha^\diamond := \alpha^\diamond_1 := Im(\alpha)$ which are to be rewritten with the help of Formulas $1(1-3)$ above. To avoid any confusion of the notation derivatives are denoted by $y' = dy/dx$ or $d^j y/dx^j$ or $y^{(j)}$, $j = 1, 2, ...$, when it is necessary. Take $\omega(x)$ satisfying the initial condition:

(8) $\omega(\alpha) = \alpha$, consequently,

(9) $p^{-1}(0) + \phi_1(\alpha^\diamond, \alpha^\diamond) = \alpha$. When $\alpha$ is fixed, the function $\phi_1$ will be written as $\phi_1(x^\diamond)$. Now we can use the change of variables Theorem 3.11 [26] for non-commutative line integrals over Cayley-Dickson algebras. Therefore,

$$(10) \quad \int_\gamma ([dy(x)/dx].h)dx = \int_\psi [dy(\omega)/d\omega].d\omega = y(\psi(b)) - y(\psi(a)),$$

where $\psi$ is a rectifiable path so that

(11) $\psi(t) = \omega(\gamma(t))$ for all $t \in [a,b]$, since by the chain rule

(12) $[dy(\omega(x))/dx].dx = [dy(\omega)/d\omega].[(d\omega(x)/dx).dx]|_{\omega=\omega(x)} = dy(\omega(x))$

(see Proposition 5 in §2). On the other hand, in view of Conditions $1(I1, I2)$ each non-commutative integration $\int_\gamma f(x)dx$ specifies a branch up to a function $\phi_f = \phi_f(\beta^\diamond{}_h, \alpha^\diamond{}_h)$ or $\phi_f(\beta^\diamond{}_h)$, when $\alpha = \gamma(a)$ is marked (fixed) and $\beta = \gamma(b)$ can vary, where

(13) $\beta^\diamond{}_h = \beta - pr_{h(\beta)}\beta$, $pr_{h(\beta)}\beta$ denotes the projection of $\beta$ on $h(\beta)$.

We have $\psi(b) = \omega(\beta)$ and $\psi(a) = \omega(\alpha)$. Suppose that $y_1(\omega(x))$ and $y_2(\omega(x) + \phi_1(x^\diamond))$ are two solutions of the differential equation (1), then their difference $y_1(\omega(\zeta)) - y_2(\omega(x) + \phi_1(x^\diamond))$ is the function of $x^\diamond$ for $Re(x) = Re(\zeta) = \alpha_0$ with $\omega(\zeta) = \omega(x) + \phi_1(x^\diamond)$, since

$[d\phi_1(x^\diamond)/dx].1 = 0$ and

$[d\{y_1(\omega(\zeta)) - y_2(\omega(x) + \phi_1(x^\diamond))\}/dx].1 =$

$[dy_1(\xi)/d\xi].([d\omega(x)/dx].1)|_{\xi=\omega(\zeta)} - [dy_2(\xi)/d\xi].([d(\omega(x)+\phi_1(x^\diamond))/dx].1)|_{\xi=\omega(x)+\phi_1(x^\diamond)}$

$= (1-1)f(\omega(\zeta)) = 0$, since the $x$-differentiability of the function $\phi_1(x^\diamond, \alpha^\diamond)$ is provided by Formulas $1(1-3)$. Thus the formula

$$(14) \quad y(\omega(x)) = \int_\alpha^x f(\omega(x))dx + \phi_f(x^\diamond)$$

gives the general solution of the (super)-differential equation (1), since the domain $U$ satisfies Conditions $1(D1, D2)$, while $\phi_f(x^\diamond) \in \mathcal{H}(U, \mathcal{A}_r)$ due to the $x$-representation $1(1-3)$, where $\omega = \omega_h$. The initial condition now is $\phi_f(\alpha^\diamond{}_h, \alpha^\diamond{}_h) = f(\alpha)$. Its partial solution is defined, when the initial condition

(15) $y(x)|_{x_0=\alpha_0} = \eta(x^\diamond)$ is given on the hyperplane

(16) $H_{\alpha_0} := \{x \in \mathcal{A}_r : x_0 = \alpha_0\}$, when $h(x)$ is not parallel to $H_{\alpha_0}$ at each point in $H_{\alpha_0}$, where

(17) $x = x_0 i_0 + x_1 i_1 + ... + x_{2^r-1} i_{2^r-1}$ is the decomposition of the Cayley-Dickson number $x \in \mathcal{A}_r$, $x_0, ..., x_{2^r-1} \in \mathbf{R}$. Here $\eta(x^\diamond)$ is a function on $\partial U$ having an $\mathcal{H}(W, \mathcal{A}_r)$ extension, where $W$ is an open neighborhood of $\partial U$.

This means that the Cauchy problem can be formulated as given by Equations $(1, 15)$ on the half-space

(18) $P_{\alpha_0} := \{x \in \mathcal{A}_r : Re(x) \geq \alpha_0\}$. Certainly, this implies $\omega(\alpha) = \alpha$ for each $\alpha \in H_{\alpha_0}$. That is $\phi_f(x^\diamond h) = \eta(x^\diamond) - \int_\alpha^x f(\tau) d\tau$ for each $x \in H_{\alpha_0}$.

Therefore, the non-commutative or non-associative analog of the Cauchy problem is given by Equations $(1, 15)$. In the previous papers [25, 25] pseudo-conformal mappings of the quaternion and octonion variables were investigated. Using the non-commutative or non-associative analog of the Riemann mapping theorem it is possible to formulate the non-commutative or non-associative analog of the Cauchy problem on more general hyper-surfaces $\partial V$ so that $\partial V = \nu(H_{\alpha_0})$, where $\nu$ is a (super-)differentiable (can be pseudo-conformal for $2 \leq r \leq 3$) diffeomorphism from $P_{\alpha_0}$ onto $V$, $V$ is the canonical closed subset in $\mathcal{A}_r$, $V = cl(Int(V))$. As usually $cl(A)$ denotes the closure of a set $A$ in $\mathcal{A}_r$, while $Int(A)$ denotes the interior of $A$ in $\mathcal{A}_r$ (see [25, 25]). For such $U = V$ the initial condition is

(19) $y(x)|_{\partial U} = \eta(x^\diamond)|_{\partial U}$, when $h(x)$ is not tangent to $\partial U$ at each point in $\partial U$, where $\eta(x^\diamond)$ is a function on $\partial U$ having an $\mathcal{H}(W, \mathcal{A}_r)$ extension, where $W$ is an open neighborhood of $\partial U$, $x^\diamond$ denotes an arbitrary point in $\partial U$. This initial condition is natural, since $y(x)$ is the function of $2^r$ real variables $x_0, ..., x_{2^r-1}$ or one Cayley-Dickson variable $x$, where $r \geq 2$.

**3. An R-linear differential equation of the first order.**

Consider a differential equation, which is $\mathbf{R}$-linear and $\mathcal{A}_r$-additive relative to an unknown function $f$ of the first order:

$$(1) \quad (dy(x)/dx).h + b(x)y(x) = Q(x)$$

on a domain $U$ satisfying Conditions $1(D1, D2)$, where a function $b(x)$ is real-valued and a mapping $Q(x)$ is octonion valued, $x \in U$, $r = 2$ or $r = 3$, $h = h(x) \in \mathcal{A}_r$. Suppose that $b(x)$, $Q(x)$ and $h(x)$ are locally analytic functions by $x$ written in the $x$-representation. When $Q$ is not identically zero this equation is called inhomogeneous. If $Q(x)$ is identically zero, then the differential equation is homogeneous:

$$(2) \quad b(x)y(x) + (dy(x)/dx).h = 0.$$

The octonion algebra is alternative, consequently,

$$(3) \quad [(dy(x)/dx).h][1/y(x)] = -b(x).$$

To solve this equation the exponential and logarithmic functions of quaternion and octonion variables can be used (see §3 in [22, 21] and §2.7 above). Since $Ln(\exp(x)) = x$, we infer

$(4)\ [dLn(\exp(x))/dx].h = h = Ih = [dLn(\xi)/d\xi].[(d\exp(x)/dx).h]|_{\xi=\exp(x)},$

where $I$ denotes the unit operator. Therefore,

$(5)\ [dLn(\xi)/d\xi]|_{\xi=\exp(x)} = [d\exp(x)/dx]^{-1}$, also $[d\exp(x)/dx].h_0 = h_0 \exp(x)$

for each $h_0 \in \mathbf{R}$, consequently,

$(6)\ [dLn(\xi)/d\xi]|_{\xi=\exp(x)}.h_0 = h_0/\xi$ and

$(7)\ [dLn\ y(\omega(x))/dx].h_0 = \{(dy(\omega(x))/d\omega).[(d\omega(x)/dx).h_0]\}[1/y(\omega(x))] = [d(Ln\ y(\omega))/d\omega].h]h_0$, since $(d\omega(x)/dx).h_0 = h(\omega)h_0$ and $(dy(\omega)/d\omega).h(\omega) = -b(z) \in \mathbf{R}$, $\omega = \omega_h(x)$. In accordance with §2 and Formula (7) Equation (3) has the general solution:

$$(8) \quad Ln\ y(\omega(x)) = -\phi_{Ln\ y}(\omega^\diamond{}_h(x)) - \int_\alpha^x b(\omega(t))dt - \phi_b(x^\diamond)$$

or

$$(9) y(\omega(x)) = \exp\{\phi_{Ln\ y}(\omega^\diamond{}_h(x)) - \int_\alpha^x b(\omega(t))dt - \phi_b(x^\diamond)\},$$

when branches of $Ln$ and $\int$ are specified.

When $b(t) \in \mathbf{R}$ for each $t$, the equality

$[d(\exp(-\int_\alpha^x b(\omega(t))dt))/dx].1 = -b(\omega(t))\exp(-\int_\alpha^x b(\omega(t))dt)$ is satisfied, since $\hat{b}(\omega).1 = b(\omega)$ and $\mathbf{R}$ is the center of the octonion algebra.

Using the variation method of the function $C = \phi_{Ln\ y}(\omega^\diamond{}_h(x))$ in (9) and $(3, 6, 7)$ one gets a solution of (1) analogously to the real case:

$$(10) \quad y(\omega(x)) = y^1(\omega(x)) + \exp\{-\int_\alpha^x b(\omega(t))dt\}$$

$$[\int_\alpha^x \exp\{\int_\alpha^\tau b(\omega(\xi))d\xi\}Q(\tau)d\tau],$$

since the octonion algebra is alternative and each equation $bx = c$ with non-zero $b$ has the unique solution $x = b^{-1}c$,

(11) $y^1(\omega(x)) = \exp\{-\int_\alpha^x b(\omega(t))dt\}\phi_1(x^\diamond)$

is any solution of the homogeneous equation, where $\omega = \omega_h(x)$ (see also Formulas 2(5 − 7)). Here $\phi_1(x^\diamond)$ is an octonion-valued function of $x^\diamond$ which can be specified by an initial condition 2(15) or more generally 2(19). When $\phi_1(x^\diamond)$ is arbitrary Formula (11) gives the general solution of the homogeneous differential equation. The second summand

$y^2(\omega(x)) := \exp\{-\int b(\omega(x))dx\}[\int \exp\{\int b(\omega(x))dx\}Q(x)dx]$

is the special solution of the inhomogeneous differential equation. Certainly, also for $r = 1$, $U \subset \mathcal{A}_1 = \mathbf{C}$, with $b(x)$ and $h(x)$ and $Q(x) \in \mathbf{C}$ for each $x \in U$ Equation (1) has the general solution (10).

Symmetrically it is possible to solve the differential equation

$$(12) \quad (dy(x)/dx).h + y(x)b(x) = Q(x)$$

on a domain $U$.

Another method of solving the inhomogeneous $\mathbf{R}$-linear equation of the first order is like to that of Bernoulli. We present $y$ as the product

(13) $y = uv$ of two super-differentiable functions $u(x)$ and $v(x)$ so that $y' = uv' + u'v$, where $u' := du(x)/dx$. With this substitution Equation (1) over quaternions becomes:

(14) $[(u'.h) + bu]v + u(v'.h) = Q(x)$. Now we put

(15) $(u'.h)(1/u) = -b$.

Its integration gives the special solution (see Formulas (3,9) above)

(16) $u(\omega(x)) = \exp\{-\int b(\omega(x))dx\}$,

where $\omega(x) = \omega_h(x)$ as in §2, branches of the indefinite integral $\int b(\omega(x))dx$ may differ on (super-)differentiable terms like $\phi(x')$ as above. It remains

(17) $u(v'.h) = Q(x)$ or $v'.h = (1/u)Q$, consequently,

(18) $v(\omega(x)) = \int \exp\{\int b(\omega(x))dx\}Q(x)dx$. Thus substituting $u$ and $v$ into the product decomposition (12) leads to:

(19) $y(\omega(x)) = y^1(\omega(x)) + \exp\{-\int b(\omega(x))dx\}\int \exp\{\int b(\omega(x))dx\}Q(x)dx$,

which is equivalent to (10), where $y^1(\omega(x))$ is given by Formula (11).

**3.1. Example.** The differential equation

(1) $[(dy(x)/dx).h]\xi(x) + y(x) = g(x)$

over octonions in the domain $Re(x) \geq \alpha_0$, where $\alpha_0 > 0$, $\xi(x)$ is any analytic real-valued function written in the $x$-representation due to Formulas $1(1-3)$, has the general solution:

(2) $\quad y(\omega_v(x)) = y^1(\omega_v(x)) + \exp(-(\omega_v(x) - \alpha))[\int_\alpha^x \exp(\omega_v(\tau) - \alpha)g(\tau)d\tau]$,

where $y^1(\omega_v(x)) = [\exp(-(\omega_v(x) - \alpha))]\phi(x^\diamond)$, while $v(x) := h(x)\xi(x)$. Particularly, one can take $\xi(x)$ in the form

$\xi(x) = \sum_k a_k x_0^{k_0}...x_{2^r-1}^{k_{2^r-1}}$,

with real expansion coefficients $a_k$, $k = (k_0, ...., k_{2^r-1})$. Further, the function $\xi(x)$ is rewritten in the $x$-representation due to Formulas $1(1-3)$.

For example, if $h(x) = x$, then the particular solution of the differential equation $2(5)$ is $\omega_h(x) = e^x$; for $h(x) = x^n$ with $n \neq 1$, $n \in \mathbf{R}$, we can take $\omega_h(x) = [(1-n)x]^{1/(1-n)}$. Though an addendum $\omega_h^1(x) = \phi_h(x^\diamond)$ can be taken so that $[d(\omega_h(x) + \omega_h^1(x))/dx].1 = h(\omega)$, since $[d\omega_h^1(x)/dx].1 = 0$.

**3.2. Note.** In components of functions Equation $3(1)$ takes the form:

(1) $\quad \sum_{j,k=0}^{2^r-1} \{i_j h_k(x) \partial y_j(x)/\partial x_k + b(x) i_j y_j(x)\} = \sum_{j=0}^{2^r-1} i_j Q_j(x)$

(see also Identities $1(1-3)$ above), where $2 \leq r \leq 3$. This is equivalent to the system of partial differential equations:

(2) $\quad \sum_{k=0}^{2^r-1} h_k(x) \partial y_j(x)/\partial x_k + b(x) y_j(x) = Q_j(x)$

for each $j = 0, 1, ..., 2^r - 1$ with real functions $y_j$, $h_j(x)$ and $b$. Initial conditions become:

(3) $\quad y_j(x)|_{\partial U} = \eta_j(x^\diamond)|_{\partial U}$ for each $j = 0, ..., 2^r - 1$. A more general system is obtained by any $\mathbf{R}$-linear transformation $(y_1, ..., y_n) = (\xi_1, ..., \xi_n)C$.

**4. A differential equation with separated variables** can be written in the form:

(1) $f(y)[dy(x)/dx].h(x) + s(x) = 0$,

where $s, h \in \mathcal{H}(U, \mathcal{A}_r)$, $f \in \mathcal{H}(W, \mathcal{A}_v)$, $2 \leq r \leq v$, $h(x)$ is a non-zero function, where $U$ and $W$ are domains satisfying Conditions $1(D1, D2)$. We seek a solution $y$ so that $y(U) \subset W$ for suitable $U$ and $W$. Integrating both sides of this equation with the help of the non-commutative line integral we deduce:

(2) $\int_\alpha^x [f(y)(dy(x)/dx).h(x)]dx + \phi_f(x^\diamond) = -\int_\alpha^x s(x)dx - \phi_s(x^\diamond)$, since

(3) $[(dy/dx).1]h(x) = (dy/dx).h(x)$

for the real-valued function $h(x)$. Due to Equations $2(5-9)$ we can make the change of variables for the integral on the right side of Equality (3).

To evaluate the integral take a function $g(x)$ satisfying Equations $1(I1, I2)$ for $f(x)$. This means that $[dg(y)/dy].1 = f(y)$, $\hat{f}(y) = dg(y)/dy$. Then $[dg(\omega(\xi))/d\xi].h = [dg(\omega)/d\omega].[(d\omega(\xi)/d\xi).h]$ for each $h \in \mathcal{A}_r$, in particular, for $\xi = g$, consequently, $(d\omega(g)/dg)^{-1} = dg(\omega)/d\omega$. One has

$[dg_j(y(x))/dx].dx = [dg_j(y)/dy].[(dy(x)/dx).dx] =$
$dg_j(y) = [(dg_j(y)/dy).dy]$,

since $d(\omega - b) = d\omega$ for each constant $b \in \mathcal{A}_r$, where $\omega(x) = \omega_h(x)$. Therefore, for the equality

(4) $\int((a_n z^n)b_n)dz = \int d(a_n z^{n+1}/(n+1))b_n$

is satisfied for each $z \in \mathcal{A}_r$, $n \neq -1$, $a_n \in \mathcal{A}_r$, $b_n \in \mathcal{A}_r$,

(5) $\int((a_n(1/z))b_n)dz = (a_n \int dLn\ z)b_n$

in $U$ not containing 0 with suitable branch of the logarithm $Ln$ and the integral $\int$ and with $a_n \in \mathbf{R}$ (see in details the integration algorithms in [21, 22, 26]). Thus with the left algorithm we infer that

(6) $\int_\alpha^x f(y)[(dy(t)/dt).h(t)]dt = \int_{t=\alpha}^{t=x} f(y)dy(\omega(t)) + \phi_f(Im\ y(\omega(x))$

for the locally analytic function $f(y)$, consequently, the equation

(7) $\int_{y(\alpha)}^{y(\omega(x))} f(y)dy = -\phi_s(x^\diamond) - \int_\alpha^x s(t)dt - \phi_f(Im\ y(\omega(x)))$

gives the relation between the function $y$ and the argument $x$ (see §2 above). When a mapping $f$ and its phrase $\nu^f$ are known, then a function $g(y) = \int f(y)dy$ and its phrase $\nu^g$ can be found easily with the help of the left algorithm. If an analytic function $f$ is with real expansion coefficients, this integration reduces to the ordinary integration. Any other branch of the integral $\int f(y)dy$ may differ from that of given by the left algorithm on a term like $\phi_f(y - Re(y)) \in \mathcal{H}$. If $f$ is a polynomial, then $g$ is a polynomial. This means that $y$ is a zero of the corresponding equation

(8) $g(y) = b(x)$,

where $b(x) = -\phi_s(x^\diamond) - \int_\alpha^x s(t)dt$.

This equation is polynomial by $y$, when $f$ is the polynomial.

Apart from the complex case zeros of polynomial equations over $\mathcal{A}_r$ generally form sets, which are unions of connected locally analytic surfaces. If polynomials are with real coefficients, real dimensions of such surfaces are up to $(2^r - 2)$. For polynomials with $\mathcal{A}_r$ coefficients which may stand on different sides of powers of $y$ dimensions of zero surfaces may be from nil up to $(2^r - 1)$ [21, 28, 37]. Thus

(9) $y = g^{-1}(b(x))$,

but generally the inverse function $g^{-1}$ need not be univalent [25] analogously to the complex case. The function $\phi_s(x^\diamond)$ may be given with the help of the boundary condition 2(19).

**5. Example.** Consider the differential equation

(1) $f(y)[dy(x)/dx].h(x) + y^m s(x) = 0$,

where $m \in \mathbf{R}$, $s, h \in \mathcal{H}(U, \mathcal{A}_r)$, $f \in \mathcal{H}(W, \mathcal{A}_r)$, $h(x)$ is real-valued, $h(x) \in \mathbf{R} \setminus \{0\}$, as in §4, $2 \leq r \leq v \leq 3$. Due to the alternativity of the octonion algebra and 4(4, 5) the equivalent equation is:

(2) $(y^{-m} f(y))[dy(x)/dx].h(x) + s(x) = 0$,

where as usually $[dy(x)/dx].h$ denotes the (super) differential of $y$ along a vector (Cayley-Dickson number) $h$. The power function $y^m$ can be written in the form $y^m = \exp(m\ Ln(y))$. Since the exponential function is analytic on $\mathcal{A}_r$ and the logarithmic function is locally analytic on $\mathcal{A}_r \setminus \{0\}$. Equation (2) is of the same type as 4(1). Thus results of §4 permit to integrate non-linear octonion (super)-differential equations. The evaluation of several useful integrals over the Cayley-Dickson algebras was given in [23, 26].

**5.1. Homogeneous differential equations.**

Let us study the homogeneous differential equation of the form:

(1) $(dy(x)/dx).h = f(yx^{-1})$ or

(2) $(dy(x)/dx).h = f(x^{-1}y)$ over octonions in the class of $x$-(super-)differentiable functions, $h = h(x) \in \mathbf{O}$. We make the substitution $y = ux$ in the first and $y = xu$ in the second case. Due to alternativity of the octonion algebra one gets from (1):

(3) $(dy/dx).h = [(du/dx).h]x + uh = f(u)$. The derivative of $u$ can be expressed from Equation (3):

(4) $(du/dx).h = (f(u) - uh)(1/x)$. Dividing on $(f(u) - uh)$ from the left we deduce:

(5) $\int [1/(f(u) - uh)][(du/dx).h] = \int (1/x)dx$. Therefore,

(6) $\phi(x^\diamond)x = \exp(\int [1/(f(u) - uh)]du(\omega_h(x)))$, consequently,

(7) $\phi(x^\diamond)x = \exp(\int [1/(f(yx^{-1}) - (yx^{-1})h)]d(yx^{-1})(\omega_h(x)))$,

where we have substituted $u = yx^{-1}$. The second differential equation can be resolved symmetrically. Thus Equality (7) presents the general integral of the first differential equation. The function $\phi(x^\diamond)$ can be rewritten as $x\phi^1(x^\diamond)$ with $\phi^1(x^\diamond) = (1/x)\phi(x^\diamond)x$, since $1/x = \tilde{x}/|x|^2$ for $x \neq 0$, $Re(x) \in \mathbf{R}$ and $\mathbf{R}$ is the center of the Cayley-Dickson algebra $\mathcal{A}_r$.

**6. Bernoulli's equation.** Let the differential equation

(1) $[dy(x)/dx].h(x) + y^k p(x) = y^m s(x)$ be given over octonions, where $h(x)$ is as in §2, $k, m \in \mathbf{R}$, $m \neq 0$, $m \neq k$. Suppose that a function $p(x)$ is real-valued, $h(x)$, $p(x)$ and $s(x)$ are $x$-differentiable on a domain $U$ satisfying Conditions $1(D1, D2)$. We change variables as $v^l = y$ for the parameter $l \in \mathbf{R} \setminus \{0\}$. When $v(x)$ and $(dv(x)/dx).h(x)$ commute for each $x \in P_\alpha$, this gives the following:

(2) $lv^{l-1}\{[dv(x)/dx].h(x)\} + v^{lk}p(x) = v^{lm}s(x)$.

The alternativity of the octonion algebra leads to:

(3) $lv^{l-lm-1}[dv(x)/dx].h(x) + v^{lk-lm}p(x) = s(x)$.

When $l = 1/(k-m)$ it reduces to

(4) $\{1/(k-m)\}v^{-1+(1-m)/(k-m)}[dv(x)/dx].h(x) + vp(x) = s(x)$.

Therefore, for $k = 1$ Equation (1) is the octonion analog of the Bernoulli's equation and due to §3 we infer:

$$(5) \quad v(\omega(x)) = v^1(\omega(x)) +$$

$$[\int_\alpha^x (1-m)s(\tau)\exp\{\int_\alpha^\tau (1-m)p(\xi)d\xi\}d\tau]\exp\{(m-1)\int_\alpha^x p(t)dt\},$$

where $\omega(x)$ is given by Equations $2(5.2, 6)$,

$v^1(\omega(x)) = \exp\{\int_\alpha^x (1-m)p(t)dt\}\phi_1(x^\diamond)$.

To avoid any misunderstanding in our notations we mention that the function $p$ from §2 is not related with the arbitrary function $p$ of this section. An arbitrary function $\phi_1 \in \mathcal{H}(U, \mathcal{A}_r)$ is determined from the boundary condition $2(15)$, so that $\phi_1(x^\diamond) = [\eta(x^\diamond)]^{1-m}$. Thus $y = v^{1/l}$ is the solution of Cauchy's problem $6(1)$, $2(15)$.

**7.** When a function $f(y)$ is $y$-differentiable on $W$, $k, m \in \mathbf{R}$, $m \neq 0$, $m \neq k$, we consider the differential equation

(1) $f(y)([dy(x)/dx].h(x)) + y^k p(x) = y^m s(x)$.

Suppose that $p(x), h(x) \in \mathbf{R} \setminus \{0\}$, $h(x)$, $p(x)$ and $s(x)$ are $x$-differentiable on $U$, open domains $U$ and $W$ in the Cayley-Dickson algebra $\mathcal{A}_r$ satisfy Conditions $1(D1, D2)$, $y(U) \subset W$. In view of Formulas $4(2-5)$ one gets the following:

(2) $\{1/(k-m)\}(v^{-1+(1-m)/(k-m)} f(v^{1/(k-m)}))([dv(x)/dx].h(x)) + vp(x) = s(x)$

for $v^l = y$ over octonions with $l = 1/(k-m)$. Thus as in §4 we deduce for $s(x) = 0$:

(3) $\beta(v) = -\phi_p(x^\diamond) - (k-m) \int_\alpha^x p(t)dt$, where

(4) $\beta(v) := \int_\eta^v z^{-2+(1-m)/(k-m)} f(z^{1/(k-m)})dz + \phi_2(Im\ v)$,

where $v = v(\omega_h(x))$, $\eta = \eta(\alpha^\diamond)$, $\alpha \in H_{\alpha_0}$. That is,

(5) $[d\beta(v(\omega_h(x)))/dx].1 = -(k-m)p(x)$.

Any other branch of the non-commutative line integral (4) may differ on a term like $\phi(v - Re(v))$. Suppose that $\beta(v)$ and its inverse function $\beta^{-1}$ are calculated. To solve (2) consider $C = -\phi_p$ as a function of $x$ so that we seek $v$ in the form:

(6) $\beta(v(\omega_h(x))) = C(x) - (k-m) \int_\alpha^x p(t)dt$.

Differentiating both sides of (6) by $x$ and using (5) we get:

$\beta^{-1}(C(x) - (k-m) \int_\alpha^x p(t)dt)[dC(x)/dx].1 = (k-m)s(x)$,

since $\beta^{-1}(\beta(v)) = v$ and $d\beta^{-1}(\xi)/d\xi = [d\beta(v)/dv]^{-1}$ for $\xi = \beta(v)$. The latter and 4(6) Formulas imply that

(7) $\psi(C(x)) = \phi_s(x^\diamond) + (k-m) \int_\alpha^x s(t)dt$, where

(8) $\psi(C(x)) := \int_{t=\alpha}^{t=x} \beta^{-1}(C(t) - (k-m) \int_\alpha^t p(\tau)d\tau)dC(t) + \phi_3(x^\diamond)$.

Thus we infer the general solution in the form:

(9) $v(\omega_h(x)) = \beta^{-1}(\psi^{-1}(\phi_s(x^\diamond) + (k-m) \int_\alpha^x s(t)dt))$, where

(10) $y = v^{1/(k-m)}$.

Arbitrary functions $\phi_3, \phi_s \in \mathcal{H}(U, \mathcal{A}_r)$ and $\phi_2 \in \mathcal{H}(W, \mathcal{A}_r)$ can be specified from the boundary condition 2(19) choosing possibly simpler formulas for $\psi$ and $\beta$:

(11) $\phi_s(\alpha^\diamond) = \psi(\beta([\eta(\alpha^\diamond)]^{(k-m)}))$.

One can mention that in the particular case of $f = 1$ and $k = 1$ the reciprocal function $\beta^{-1} = \exp$ is the exponential function and Equations $(9, 10)$ coincide with $3(10, 11)$ for $h = 1$.

**8. The total differential equations.**

If $z = F(x, y)$ is a (super-)differentiable function of two Cayley-Dickson variables $x, y$, then

(1) $dz = [\partial F(x, y)/\partial x].dx + [\partial F(x, y)/\partial y].dy$, where $\partial F(x, y)/\partial x$ and $\partial F(x, y)/\partial y$ are **R**-homogeneous $\mathcal{A}_r$-additive operators. At the same time each super-differentiable function is infinite differentiable, consequently,

(2) $[\partial^2 F(x, y)/\partial x \partial y].(h, v) = [\partial^2 F(x, y)/\partial y \partial x].(v, h)$ for all Cayley-Dickson numbers $h, v \in \mathcal{A}_r$.

The family of all **R**-homogeneous $\mathcal{A}_r$-additive operators on $\mathcal{A}_r^m$ with values in $\mathcal{A}_r^n$ is denoted by $L_q(\mathcal{A}_r^m, \mathcal{A}_r^n)$ for natural numbers $m, n \in \mathbf{N}$, $2 \leq r$. Suppose that there is given the differential 1-form over the Cayley-Dickson algebra $\mathcal{A}_r$:

(3) $w = A(x, y).dx + B(x, y).dy$,

where $A, B \in \mathcal{H}(U \times V, L_q(\mathcal{A}_r^2, \mathcal{A}_r))$, $U$ and $V$ are open domains in $\mathcal{A}_r$ satisfying Conditions $1(D1, D2)$. In accordance with $(1, 2)$ it is the complete differential if and only if

(4) $[\partial A(x, y)/\partial y].(h, v) := [\partial(A(x, y).h)/\partial y].v$
$= [\partial(B(x, y).v)/\partial x].h =: [\partial B(x, y)/\partial x].(h, v)$ for all $h, v \in \mathcal{A}_r$.

Mention that $d^2 z = \sum_{j=0}^{2^r-1} i_j d^2 z_j = 0$ for the Cayley-Dickson variable $z$. The condition $dF = 0$ is equivalent to $F(x, y) = C$, where $C$ is a Cayley-Dickson constant. In this case $dy = -(\partial F/\partial y)^{-1}.[(\partial F/\partial x).dx]$, where the composition of two $L_q$ operators $G \in L_q(\mathcal{A}_r^n, \mathcal{A}_r^m)$, $K \in L_q(\mathcal{A}_r^p, \mathcal{A}_r^n)$ is defined as usually $G.[K.h] = G(K(h))$ for mappings for each $h \in \mathcal{A}_r^p$, $p, n, m \in \mathbf{N}$. If Condition (4) is satisfied, we seek a solution $f$ of (3) as $w = dF$. This means that $\partial F/\partial x = A$ so that $F(x, y) = \int_\alpha^x A(t, y).dt + G(x^\diamond, y)$. From the condition $\partial F/\partial y = B$ it follows that

$B(x, y).h = \int_\alpha^x (\partial A(t, y)/\partial y).(dt, h) + (\partial G(x^\diamond, y)/\partial y).h$. Using Condition (4) we rewrite the latter equation as

$B(x, y).h = \int_\alpha^x (\partial B(t, y)/\partial t).(dt, h) + (\partial G(x^\diamond, y)/\partial y).h$. After the integration one gets:

$B(x, y) - B(\alpha, y) + \partial G(x^\diamond, y)/\partial y = B(x, y)$, consequently, $\partial G(x^\diamond, y)/\partial y =$

$B(\alpha, y)$ and inevitably

$G(x^\diamond, y) = \int_\beta^y B(\alpha, \tau).d\tau + \phi_B(x^\diamond, y^\diamond)$. Thus

$$(5) \quad F(x,y) = \int_\alpha^x A(t,y).dt + \int_\beta^y B(\alpha,\tau).d\tau + \phi_B(x^\diamond, y^\diamond),$$

where $\alpha, x \in U$, $\beta, y \in V$. Therefore, Condition (4) is sufficient for the existence of the function $F$. A boundary condition can be specified as:

(6) $F(x,y)|_{x_0=\alpha_0, y_0=\beta_0} = \psi(x^\diamond, y^\diamond)$.

Initially domains $U$ and $V$ can be taken as half-spaces $\{z : Re(z) \geq \tau_0\}$ with either $\tau_0 = \alpha_0$ or $\tau_0 = \beta_0$ respectively and then they can be transformed with the help of pseudo-conformal or more generally $z$-(super-)differentiable diffeomorphisms into given domains in accordance with the octonion analog of the Riemann mapping theorem (see also §2 above).

**8.1. Remark.** Generally in §§4-7 differential equations can be non-linear. In real coordinates 8(3) becomes the system of real differential forms:

(1) $w_j = \sum_{k=0}^{2^r-1} \{a_{j,k}(x,y)dx_k + b_{j,k}(x,y)dy_k\}$ for every $j = 0, 1, ..., 2^r - 1$,

where $A(x,y).i_k = \sum_{j=0}^{2^r-1} i_j a_{j,k}(x,y)$ and $B(x,y).i_k = \sum_{j=0}^{2^r-1} i_j b_{j,k}(x,y)$ with real-valued functions $a_{j,k}(x,y)$ and $b_{j,k}(x,y)$. For them Condition 8(4) reads as:

(2) $\partial a_{j,k}(x,y)/\partial y_l = \partial b_{j,l}(x,y)/\partial x_k$ for all $j, k, l = 0, 1, ..., 2^r - 1$.

**8.2. Example.** Let $A(x,y).dx = (x^2 + xy + yx)dx + (dx)(x^2 + xy + yx) + (x(dx)x + x(dx)y + y(dx)x)$, $B(x,y).dy = (x^2 - y^2)dy + (dy)(x^2 - y^2) + (x(dy)x - y(dy)y)$, then $\partial A(x,y)/\partial y = (xI_2 + I_2 x)I_1 + I_1(xI_2 + I_2 x) + (xI_1 I_2 + I_2 I_1 x) = \partial B(x,y)/\partial x$, where $I$ denotes the unit operator, $(I_1 I_2).(h,v) := hv$, $(I_2 I_1).(h,v) := vh$ for all $h, v \in \mathcal{A}_r$. The line integration gives the function $F$:

$$(2) \quad F(x,y) = \int_\alpha^x \{(t^2+ty+yt)dt+(dt)(t^2+ty+yt)+(t(dt)t+t(dt)y+y(dt)t)\}$$

$$+ \int_\beta^y \{(\alpha^2 - \tau^2)d\tau + (d\tau)(\alpha^2 - \tau^2) + (\alpha(d\tau)\alpha - \tau(d\tau)\tau)\} + \kappa(x^\diamond, y^\diamond) =$$

$$x^3 + x^2 y + yx^2 + xyx - \alpha^3 - \alpha^2 y - y\alpha^2 - \alpha y\alpha$$

$$+\alpha^2 y + y\alpha^2 + \alpha y\alpha - y^3 - \alpha^2 \beta - \beta\alpha^2 - \alpha\beta\alpha + \beta^3 + \kappa(x^\diamond, y^\diamond)$$

$= (x^3 + x^2 y + yx^2 + xyx - y^3) + \phi(x^\diamond, y^\diamond)$, where the numbers $\alpha$ and $\beta$ are marked, $\alpha \in \partial U$, $\beta \in \partial V$, while $x, y$ are variable Cayley-Dickson numbers.

**8.3. Example.** Suppose that

$$A(x,y) = (\sum_{n=0}^{\infty}[(-1)^n/(2n+1)!]\sum_{k=0}^{2n} x^k I_1 x^{2n-k})\cos(y),$$

$$B(x,y) = \sin(x)(\sum_{n=1}^{\infty}[(-1)^n/(2n)!]\sum_{k=0}^{2n-1} y^k I_2 y^{2n-1-k}).$$

Evidently, Condition 8(4) is accomplished. From Formula 8(5) one deduces:

$$F(x,y) = \sin(x)\cos(y) + \phi(x^\diamond, y^\diamond).$$

**9. Integration factor.** Only in rare cases the differential $A(x,y).dx + B(x,y).dy$ is complete. We consider the problem of an integration factor existence so that

(1) $w = \mu A.dx + \mu B.dy$ would be the complete differential, where $\mu \in \mathcal{H}(U \times V, \mathcal{A}_r)$, $2 \le r \le 3$. This means that

(2) $(\partial(\mu(x,y)A(x,y))/\partial y).(h,v) = (\partial(\mu(x,y)B(x,y))/\partial x).(h,v)$ for all $h, v \in \mathcal{A}_r$. That is

(3) $\mu(x,y)([\partial(A(x,y).h)/\partial y].v - [\partial(B(x,y).v)/\partial x].h) = -([\partial\mu(x,y)/\partial y].v)(A(x,y).h) + ([\partial\mu(x,y)/\partial x].h)(B(x,y).v)$

for all $h, v \in \mathbf{O}$ which can be $h = i_j$ and $v = i_k$ with $j, k = 0, 1, ..., 2^r - 1$. This is the system of partial differential equations of the first order for an unknown integration factor $\mu(x,y)$. Now we consider particular cases of (3).

Suppose that $\mu$ depends only on $x$, hence $\partial\mu/\partial y = 0$. Therefore, $\mu$ satisfies the identity:

(4) $(1/\mu(x))\{([\partial\mu(x)/\partial x].h)(B(x,y).v)\} = ([\partial(A(x,y).h)/\partial y].v - [\partial(B(x,y).v)/\partial x].h)$

over octonions for all $h = i_j$ and $v = i_k$ with $j, k = 0, 1, ..., 2^r - 1$, when $\mu(x) \ne 0$. The latter equality simplifies over quaternions due to the associativity law as:

(5) $(1/\mu(x))([\partial\mu(x)/\partial x].h) = ([\partial(A(x,y).h)/\partial y].v - [\partial(B(x,y).v)/\partial x].h)(B(x,y).v)^{-1}$.

On the left of Equation (5) stands the function of $x$ only, consequently, the right side also may depend on $x$ only. Its integration for each $h = i_j$ and the summation by $j = 0, ..., 3$ gives the solution:

(6) $\mu(x) = \exp\{\sum_{j=0}^{3}\int_{(x_0,...,x_{j-1},x_{j,0},x_{j+1,0},...,x_{3,0})}^{(x_0,...,x_j,x_{j+1,0},...,x_{3,0})} ([\partial(A(x,y).i_j)/\partial y].v$

$-[\partial(B(x,y).v)/\partial x].i_j)(B(x,y).v)^{-1}dx_j\},$

where $v = i_k$, $k$ is a constant or $k = k(j) \in \{0, 1, 2, 3\}$, $x = x_0 i_0 + ... + x_3 i_3 \in \mathbf{H}$, $x_j \in \mathbf{R}$ for each $j$, $(x_{0,0}, ..., x_{3,0})$ is a marked point. For example, if $([\partial(A(x,y).i_j)/\partial y].v - [\partial(B(x,y).v)/\partial x].i_j)(B(x,y).v)^{-1} = i_j$ for each $j$, then $\mu = e^x$.

For $\mu$ depending only on $y$ Equation (3) simplifies:

(7) $\mu(y)([\partial(A(x,y).h)/\partial y].v - [\partial(B(x,y).v)/\partial x].h) = -([\partial\mu(y)/\partial y].v)(A(x,y).h)$.

Particularly, over quaternions it can be written as:

(8) $[1/\mu(y)]([\partial\mu(y)/\partial y].v) = (-[\partial(A(x,y).h)/\partial y].v + [\partial(B(x,y).v)/\partial x].h)(A(x,y).h)^{-1}$

for $\mu(y) \neq 0$. Taking $v = i_k$ and integrating by $y_k$ we infer that:

(9) $\mu(y) = \exp\{\sum_{k=0}^{3} \int_{(y_0,...,y_{k-1},y_{k,0},y_{k+1,0},...,y_{3,0})}^{(y_0,...,y_k,y_{k+1,0},...,y_{3,0})} ([\partial(B(x,y).i_k)/\partial x].h$
$- [\partial(A(x,y).h)/\partial y].i_k)(A(x,y).h)^{-1} dy_k\}$,

where $h = i_j$, $j$ is a constant or $j = j(k) \in \{0, 1, 2, 3\}$.

If $\mu(x,y) = \mu(xy)$ is the function of $z = xy$, then $[\partial\mu(xy)/\partial x].h = [d\mu(z)/dz].(hy)$ and $[\partial\mu(xy)/\partial y].h = [d\mu(z)/dz].(xh)$ over octonions. In this case Equation (3) becomes again simpler:

(10) $\mu(xy)([\partial(A(x,y).h)/\partial y].v - [\partial(B(x,y).v)/\partial x].h) = -([d\mu(z)/dz].(xv)(A(x,y).h) + ([d\mu(z)/dz].(hy)(B(x,y).v)$

for all $h, v \in \mathbf{O}$, particularly, for $h = i_j$ and $v = i_k$ with $j, k = 0, 1, ..., 2^r - 1$. Due to the alternativity of the octonion algebra the latter identity can be transformed into:

(11) $[1/\mu(xy)]\{([d\mu(z)/dz].(hy)(B(x,y).v) - ([d\mu(z)/dz].(xv)(A(x,y).h)\} = ([\partial(A(x,y).h)/\partial y].v - [\partial(B(x,y).v)/\partial x].h)$,

when $\mu(xy) \neq 0$, where $z = xy$. For any marked non-zero $x$ and the variable $y$ (or vice versa) Equation (11) can be integrated by the variable $z$ as in §4, since $y = x^{-1}z$ in the octonion algebra, when $x \neq 0$.

Certainly, over quaternions and octonions more general integration factors $\mu$ and $\nu$ can be considered:

(11) $dF = \mu[(A.dx)\nu] + \mu[(B.dy)\nu]$ because this is equivalent to the equality:

(12) $[(1/\mu)dF](1/\nu) = (A.dx) + (B.dy)$,

when $\mu \neq 0$ and $\nu \neq 0$. But generally the differential form $w = \{\sum_{j=0}^{n} \mu_j[(A.dx)\nu_j] + \mu_j[(B.dy)\nu_j]\}$ may not be equivalent to $A.dx + B.dy$.

**10. Differential equations of the $n$-th order by $dy/dx$.** Let

(1) $\sum_{k=0}^{n}\{A_k, ([dy(x)/dx].h)\} = 0$

be the differential equation of the $n$-th order by $dy/dx$, where $h = h(x, y) \in \mathcal{A}_r$, $2 \le r$, $\{A_k, z\} := \{a_{k,1} z a_{k,2} z ... z a_{k,k+1}\}_{q(2k+1)}$ for $z \in \mathcal{A}_r$, $a_{k,j} = a_{k,j}(x, y)$ are coefficients, $q(2k+1)$ are vectors indicating on an order of multiplications. For $r = 2$ due to the associativity of the quaternion skew field vectors $q(2k+1)$ can be omitted. Suppose that some solution $[dy(x)/dx].h = {}_\alpha f(x, y)$ of (1) is found. Than the latter equation is the first order differential equation relative to the derivative $dy/dx$. Apart from the complex case over the Cayley-Dickson algebra the decomposition $\sum_{k=0}^{n}\{A_k, z\} = (z - {}_1 f(x, y))...(z - {}_n f(x, y))$ may not be existing. If it is the case, then the initial differential equation reduces to $n$ differential equations of the first order relative to the derivative $dy/dx$: $[dy(x)/dx].h = {}_l f(x, y)$, $l = 1, ..., n$. One such equation over the Cayley-Dickson algebra $\mathcal{A}_r$ gives the corresponding system of $2^r$ differential equations over the real field $\mathbf{R}$.

**11. Example.** The differential equation

(1) $x([dy(x)/dx].h)^2 - y(x)([dy(x)/dx].h) - x\{([dy(x)/dx].h)x^{-1} y(x)\} + x = 0$,

where $h = h(x)$, over quaternions after the multiplication on $x^{-1}$ from the left takes the form

(2) $([dy(x)/dx].h)^2 - (x^{-1} y(x))([dy(x)/dx].h) - ([dy(x)/dx].h)(x^{-1} y(x)) + 1 = 0$,

where $h = h(x) \in \mathbf{H}$. The quadratic polynomial $z^2 + bz + zb + c$ over the Cayley-Dickson algebra $\mathcal{A}_r$ can be presented as $(z - \lambda_1)(z - \lambda_2)$, when $bz + zb = -\lambda_1 z - z\lambda_2$ and $\lambda_1 \lambda_2 = c$. Thus

(3) $z^2 + bz + zb + c = (z + b)^2 + (c - b^2) = 0$

and the zeros of the quadratic polynomial are the following:

(4) $\lambda_1 = -b - (b^2 - c)^{1/2}$, $\lambda_2 = -b + (b^2 - c)^{1/2}$.

Equation (2) can be considered over octonions as well. Applying the latter formula to (2) one gets:

(5) $[dy(x)/dx].h = x^{-1} y(x) \pm \{(x^{-1} y(x))^2 - 1\}^{1/2}$.

Each Cayley-Dickson non-zero number $z \in \mathcal{A}_r$, $2 \le r$, can be written in the polar form $z = |z|e^{M\phi}$, where $\phi \in \mathbf{R}$, $M \in \mathcal{A}_r$ is a purely imaginary Cayley-Dickson number with $|M| = 1$ (see §3 in [21, 22]). Therefore, $z^{1/2} =$

$\rho e^{K\psi}$, where $\rho > 0$, $K \in \mathcal{A}_r$ is a purely imaginary Cayley-Dickson number, $|K| = 1$, $\psi \in \mathbf{R}$, so that $\rho^2 = |z|$, $e^{2K\psi} = e^{M\phi}$ or $\cos(2\psi) + K\sin(2\psi) = \cos(\phi) + M\sin(\phi)$. That is, $\cos(2\psi) = \cos(\phi)$, $K\sin(2\psi) = M\sin(\phi)$. For $M\sin(\phi) \neq 0$ it gives two solutions: $K = M$, $\psi = \phi/2$ $(mod\ 2\pi)$ or $\psi = \phi/2 + \pi$ $(mod\ 2\pi)$. For $M\sin(\phi) = 0$, i.e. $\phi = 0$ $(mod\ \pi)$, the Cayley-Dickson number $K$ may be arbitrary purely imaginary of the unit norm, $\psi = \phi/2$ $(mod\ \pi)$. For $\phi = 0$ $(mod\ 2\pi)$ one gets two solutions $\rho$ and $-\rho$, since $e^{K\pi} = -1$. For $\phi = \pi$ $(mod\ 2\pi)$ this gives the imaginary unit sphere of solutions, since $e^{K\pi/2} = K$, $e^{3K\pi/2} = -K$, $K^2 = -1$.

**12. Lagrange's and Clairaut's differential equations.**

Let us consider the differential equation

(1) $y = f(x, (dy(x)/dx).h)$, where $f(x, p)$ is a (super-)differentiable function of $x$ and $p$ on open domains $U$ and $V$ in the Cayley-Dickson algebra $\mathcal{A}_r$, $2 \leq r$, $x \in U$, $p \in V$, $h = h(x)$ or $h = h(x, y) \in \mathcal{A}_r$ is a constant or a (super-)differentiable function, $x \in U$, $y \in W$. We put $p := (dy(x)/dx).h$. After the (super-)differentiation by $x$ one obtains:

(2) $p = [\partial f(x,p)/dx].h + [\partial f(x,p)/\partial p].[(dp/dx).h]$ or

(3) $[\partial f(x,p)/\partial p].[(dp/dx).h] + \{[\partial f(x,p)/\partial x].h - p\} = 0$. Now we have the differential equation of the first order by $p$ as the function by $x$. Suppose that an integral

(4) $\Phi(x, p, C) = 0$ of Equation (3) is found, where $C = C(\omega^\diamond{}_h)$ may be a (super-)differentiable function, the mapping $\omega^\diamond{}_h = Im\ \omega_h$ is as above (see Equations $2(5-6)$). Then it is possible to eliminate the parameter $p$. From the equality $y = f(x, p)$ we infer:

(5) $(dp(x)/dx).h = (\partial f/\partial p)^{-1}.\{p - [\partial f/\partial x].h\}$,

where $(\partial f/\partial p)^{-1}$ denotes the inverse operator whenever it exists. The given relations $(4, 5)$ between $x, y$ and $C$ provide the general integral of (1).

As Lagrange's equation we undermine the differential relation between $x$ and $y$ which is $\mathbf{R}$-linear and $\mathcal{A}_r$-additive by $x$ and $y$:

(6) $y = xf(p) + s(p)x + \eta(p)$,

where $p = [dy(x)/dx].h$, $h = h(x)$ or $h = h(x,y)$, $f$ and $s$ and $\eta$ are super-differentiable functions on suitable open domains in the Cayley-Dickson algebra. We consider the case when only one of the functions either $f(p)$ or $s(p)$ are not identically zero. After the (super-)differentiation of the latter equality by $x$ we deduce that

(7) $p = hf(p) + s(p)h + \{x((df(p)/dp).[(dp(x)/dx).h]) + ((ds(p)/dp).[(dp(x)/dx).h])x\}$
$+ [d\eta(p)/dp].[(dp(x)/dx).h]$
$= hf(p) + s(p)h + \{x(df(p)/dp).() + (ds(p)/dp).()x + [d\eta(p)/dp].()\}.[(dp(x)/dx).h].$

Expressing $dx/dp = (dp/dx)^{-1}$ from (7) we obtain the equation:

(8) $(dp/dx).h = \{x(df(p)/dp).() + (ds(p)/dp).()x + [d\eta(p)/dp].()\}^{-1}([p - ()f(p) - s(p)()].h)$ or

(9) $(dx/dp).h = [p - ()f(p) - s(p)()]^{-1}[\{x(df(p)/dp).() + (ds(p)/dp).()x + [d\eta(p)/dp].()\}.h],$

where $h$ on the right is substituted into the brackets () of the composition of operators whenever the inverse operator $[p - ()f(p) - s(p)()]^{-1}$ exists. One treats here $p$ as the independent variable and $x$ is the dependent variable. The latter differential equation is of the first order by $x$. It may be solved by the method either of Bernoulli or of variable constants as above in §§3 and 6. The general solution should be $x = \Phi(p, C)$, where $C = C(Im\ p_h)$ is a super-differentiable function. Substituting this expression into Lagrange's equation gives $y$ as the function of $p$ and $C(Im\ p_h)$.

Any relation of the form $x = \Phi(p, C)$ or $y = \Psi(p, C)$ supplies parametrix by $C(Im\ p_h)$ representation of a seeked integral. The auxiliary parameter $p$ should be excluded.

Clairaut's differential equation over octonions is the following:

(10) $y = x[(dy/dx).h] + \eta((dy/dx).h) = xp + \eta(p)$

consists in this particular case of the function $f(p) = p$, where $h = 1$. Therefore, after the differentiation we deduce that

(11) $(1 - h)p = (d\eta/dp).[(dp/dx).h] + x[(dp/dx).h]$, consequently,

(12) $[(d\eta/dp) + xI].dp = 0,$

since $h = 1$, where the $I$ denotes the unit operator. The latter equation can be satisfied, if $(dp/dx).1 = 0$, hence $p(x) = \phi(x^\diamond)$. The substitution of $p$ in the initial Clairaut's equation gives:

(13) $y(x) = x\phi(x^\diamond) + \eta(\phi(x^\diamond)).$

The second possibility provided by Equality (12) is: $(xI + d\eta/dp) = 0$, that leads to the solution in the parameter form:

(14) $xI = -d\eta/dp$ and $y = -(d\eta/dp).p + \eta(p)$,

where $p$ plays the role of the parameter.

For $h \neq 1$ in Clairaut's differential equation one can use at first the substitution of the variable $x$ onto $t$ reducing it to the considered case of $h = 1$, i.e. $(d\omega_h(t)/dt).1 = h(t)$, where the function $\omega_h(t)$ was considered above (see Formula 2(5)). Then $(dy/dt).1 = (dy/d\omega_h).[(d\omega_h(t)/dt).1]$ for $x = \omega_h(t)$. Making this substitution one should either fix a representation for $\omega_h(t)$ to avoid ambiguity or consider the general solution.

**12.1. Example.** Let us consider Lagrange's differential equation over octonions:

(1) $y(x) = ((dy(x)/dx).h)(2x+1)$.

After the differentiation by $x$ one gets:

(2) $p(1-2h) = [(dp/dx).h](2x+1)$, where $p = (dy/dx).h$, consequently,

(3) $[2x+1]^{-1}[(dx/dp).h] = (1-2h)^{-1}(1/p)$, when $h \neq 1/2$. If $h = const \neq 1/2$ is the marked octonion, then

(4) $y(\omega_h(x)) = \int_\alpha^x \exp\{[(1-2h)/2]Ln(2\omega_h(t)+1)\}dt + \phi_y(Im\ \omega_h(x))$.

If $(1-2h) = s(x)(2x+1)$, then

(5) $y(\omega_h(x)) = \int_\alpha^x \exp\{s(\omega_h(t))\}dt + \phi_y(Im\ \omega_h(x))$.

**12.2. Example.** Consider Clairaut's differential equation over octonions:

(1) $y = x[(dy/dx).1] - [(dy/dx).1]^2/4$.

In accordance with Equations 12(13, 14) we obtain two solutions:

(2) $y = x\phi(x^\diamond) - [\phi(x^\diamond)]^2/4$ and

(3) $y = p^2/2 - p^2/4 = p^2/4$ and $xp = p^2/2$, so that $x = p/2$, consequently, $y = x^2$. The second solution is correct, since $p$ and $x$ commute, consequently, also $p$ and $(dp/dx).1$ commute.

**12.3. Example.** The next equation that we consider is Lagrange's differential equation over octonions:

(1) $y = (x+1)((dy/dx).1)^2$.

With the notation $(dy/dx).1 = p$ this equation is equivalent to $y = xp^2 + p^2$. We shall seek a solution in the case, when $p = (dy/dx).1$ and $(dp/dx).1$ commute with each other. This leads to the commutativity of $p$ with $(dp/dz).1$, consequently,

(2) $(p^2 - p) + 2p(x+1)(dp/dx).1 = 0$

in accordance with §12. Then one reduces this equation to the linear relative to $x$:

(3) $(dx/dp).1 + (2/(p-1))x = 2/(1-p)$,

when $p^2 - p \neq 0$. Its integration gives:

(4) $x = [\phi(p - Re(p))/(p-1)^2] - 1$

and its substitution into the first equation produces

(5) $y = \phi(p - Re(p))p^2/(p-1)^2$, when $\phi(p - Re(p))$ and $p$ commute,

for example, when $\phi$ is real-valued or $Im[\phi(p - Re(p))]$ belongs to $Im(p)\mathbf{R}$, where $Im(p) := p - Re(p)$ denotes the imaginary part of $p$. Thus $(4, 5)$ is the solution in the parameter form. Excluding the parameter $p$ from these two equations we obtain the solution

(6) $y = [(x+1)^{1/2} + C]^2$, where $C^2 = \phi(Im(p))$.

It remains to consider the equation $p^2 - p = 0$. It has two solutions $p = 0$ and $p = 1$, which lead to two solutions of Equation (1): the singular solution $y = 0$ and the particular solution $y = x + 1$.

Let also

(7) $y = (1 + ((dy/dx).1))x + ((dy/dx).1)^2$

be Lagrange's equation over octonions. Seeking its solution for commuting $p$ and $x$ with the help of §12 one infers:

(8) $x = \phi(Im(p))e^{-p} - 2p + 2$ and $y = \phi(Im(p))(1+p)e^{-p} - p^2 + 2$, when $Im[\phi(Im(p))]$ belongs to the straight line $Im(p)\mathbf{R}$.

**13. Differential equations of the second and higher order.**

A general differential equation of the $n$-th order over the Cayley-Dickson algebra $\mathcal{A}_r$ can be written in the form:

(1) $F(x, y, y^{(1)}.h_{1,1}, ..., y^{(n)}.(h_{n,1}, ..., h_{n,n})) = 0$,

where $F$ is a function, $x \in U$, an open domain $U$ satisfies conditions of §1. We consider $\mathbf{R}$-linear $\mathcal{A}_r$-additive differential equations of the particular expression:

(2) $a_n(x)y^{(n)}.(h_{n,1}, ..., h_{n,n}) + a_{n-1}(x)y^{(n-1)}.(h_{n-1,1}, ..., h_{n-1,n-1}) + ... + a_1(x)y^\diamond.h_{1,1} + a_0(x)y = g(x)$,

where $h_{j,k} = h_{j,k}(x)$, $1 \leq k \leq j \leq n$, $h_{j,k}, a_0, ..., a_n, g$ are (super-)differentiable functions, $y' = y^{(1)} = dy/dx, ..., y^{(n)} = d^n y/dx^n$, $y^{(0)} = y$. Special cases of this differential equation are considered below permitting to obtain a solution. Let

(3) $(...(y^{(n)}.h_1)...).h_n = g(x)$,

where $(...(y^{(n)}.h_1)...).h_n = [d(...(y^{(n-1)}.h_1)...).h_{n-1}/dx].h_n$.

The first integration gives the relation:

(4) $(...(y^{(n-1)}(\omega_{h_n}(x)).h_1)...).h_{n-1} = -\phi_{y,1}(Im\ \omega_{h_n}) + \psi_1(x)$,

where $\psi_1(x) = \phi_g(x^\diamond) + \int_\alpha^x g(t)dt$, $\alpha$ and $x \in U$ (see also §2), $\phi_g(x^\diamond) = \phi_g(x^\diamond, \alpha^\diamond)$, $\phi_{y,1}(Im\ \omega_{h_n}) = \phi_{y,1}(Im\ \omega_{h_n}, \alpha^\diamond{}_{h_n})$, $\psi_1(x) = \psi_1(x, \alpha)$, but parameter $\alpha$ is omitted for short for some marked Cayley-Dickson number $\alpha$. It is convenient to use here the left algorithm for calculation of the line integral over the Cayley-Dickson algebra. The second integration leads to:

(5) $y^{(n-2)}(\omega_{h_{n-1}}(\omega_{h_n}(x))).(h_1,...,h_{n-2}) = -\phi_{y,2}(Im\ \omega_{h_{n-1}}) + x\{-\phi_{y,1}(Im\ \omega_{h_n}) + \phi_g(x^\diamond)\} + \psi_2(x)$, where $\psi_2(x) = \int_\alpha^x (\psi_1(t) - \phi_g(t^\diamond))dt$.

Repeating such non-commutative integration $n$ times one gets the result:

(6) $y(\omega_{h_1}(...(\omega_{h_n}(x)))) = \psi_n(x) + x^{n-1}\{-\phi_{y,1}(Im\ \omega_{h_n}) + \phi_g(x^\diamond)\}/(n-1)! + x^{n-2}\{-\phi_{y,2}(Im\ \omega_{h_{n-1}}) + \phi_{\psi_1}(x^\diamond)\}/(n-2)! + ... + \{-\phi_{y,n}(Im\ \omega_{h_1}) + \phi_{\psi_{n-1}}(x^\diamond)\}$,

where $\psi_k(x)$ is given by the recurrence relation:

(7) $\psi_k(x) = \int_\alpha^x (\psi_{k-1}(t) - \phi_{\psi_{k-1}}(t^\diamond))dt + \psi_k(x^\diamond)$ for each $k = 2,...,n$.

Functions $\{-\phi_{y,1}(Im\ \omega_{h_n}) + \phi_g(x^\diamond)\}$ and $\{-\phi_{y,k}(Im\ \omega_{h_{n-1}}) + \phi_{\psi_{k-1}}(x^\diamond)\}$ for $k = 2,...,n$ can be determined from the initial Cauchy conditions

(8) $[(...(y^{(k)}(x).h_1)...).h_k]|_{\partial U} = \eta_k(x^\diamond)|_{\partial U}$ for each $k = 0,1,...,n-1$,

when $h_j(x)$ is not tangent to $\partial U$ at each point $x^\diamond$ in $\partial U$ for each $j = 1,...,n$. Here the $\eta_k(x^\diamond)$ is a function on $\partial U$ having an $\mathcal{H}(W, \mathcal{A}_r)$ extension, where $W$ is an open neighborhood of $\partial U$. With the initial conditions this gives the solution in the domain $U = P_{\alpha_0}$:

(9) $y(\omega_{h_1}(...(\omega_{h_n}(x)))) = \eta_0(_0x) + [(x - {}_0x)\eta_1(_0x) + ... + (x - {}_0x)^{n-1}/(n-1)!]\eta_{n-1}(_0x) + \int_{0x}^x (x-t)^{n-1}g(t)dt/(n-1)!$,

where $P_{\alpha_0}$ is given by Condition 2(18), ${}_0x = \alpha_0 + x - Re(x)$ is the point in $H_{\alpha_0} = \partial P_{\alpha_0}$.

When the following equation

(10) $F(x, y^{(1)}.h_1, ..., (...(y^{(n)}.h_1)...).h_n) = 0$

is given not depending explicitly on $y$ with (super-)differentiable functions $h_k = h_k(x)$ for each $k = 1,...n$, then one can replace $y^{(1)}.h_1$ with a function

(11) $v = y^{(1)}.h_1$. The (super-)derivative of $v(x)$ is

(12) $[dv(x)/dx].h_2 = \{d[(dy(x)/dx).h_1]/dx\}.h_2 = [d^2y(x)/dx^2].(h_1, h_2) + (dy(x)/dx).[(dh_1(x)/dx).h_2]$

in accordance with Proposition 2.6. There exists $p_1(x)$ so that $h_1(x) = p_1(x) + [dp_1(x)/dx].h_2$ as the solution of such first order differential equation (see §3). Thus by induction we get the following:

(13) $(...(y^{(k+1)}.h_1)...).h_{k+1} = (...(v^{(k)}.h_2)...).h_{k+1}$ for each $k = 1, ..., n-1$.

With the new function $v$ Equation (10) becomes

(14) $F(x, v, v^{(1)}.h_2, ..., (...(v^{(n-1)}.h_2)...).h_n) = 0$.

If any solution $v$ is obtained, then

(15) $y(\omega_{h_1}(x)) = \int v(x)dx$.

Practically $(1, 2, 3)$ and $(10)$ in real coordinates correspond to the specific systems of $2^r$ partial differential equations of $2^r$ real variables with generally non-constant coefficients, because $a_k$ and $h_{j,k}$ and $h_j$ may be non-constant functions.

For the differential equation of the form

(16) $F(y, y^{(1)}.h_1, ..., (...(y^{(n)}.h_1)...).h_n) = 0$

with super-differentiable functions $h_k = h_k(x)$ for each $1 \le k \le n$ we express $y^{(1)}$ as a function of $y = y(x)$. For this we put $u = y^{(1)}.h_1$. The (super-)derivative of this function is:

(17) $[du/dx].h_2 = [du/dy].[(dy/dx).h_2] = [du/dy].s_1 = [d^2y(x)/dx^2].(h_1, h_2) + (dy(x)/dx).[dh_1(x)/dx].h_2$,

where $s_1 = s_1(x) = (dy(x)/dx).h_2$. There exists a solution $p_1$ satisfying the differential equation $h_1(x) = p_1 + [dp_1(x)/dx].h_2(x)$ (see §3). Continuing this recurrence relations by induction one obtains the new differential equation relative to $u$:

(18) $F(u, u^{(1)}.s_1, ..., (...(u^{(n-1)}.s_1)...).s_{n-1}) = 0$,

where $s_j = s_j(x) = (dy(x)/dx).h_{j+1}$ for each $j$.

If a general solution $u$ of (18) is calculated, then

(19) $\int u^{-1}([dy/dx].h_1)dx = x + \phi_1(x^\diamond)$,

since $u^{-1}([dy(x)/dx].h_1) = 1$ for each $x$.

For the differential equation

(20) $(...(y^{(n)}.h_1)...).h_n = g((...(y^{(n-1)}.h_1)...).h_{n-1})$,

where $h_j = h_j(x)$ for each $j$, it is useful to put

(21) $u(x) = (...(y^{(n-1)}.h_1)...).h_{n-1}$.

For this function its (super-)derivative is given by the equality:

(22) $(du(x)/dx).h_n = (...(y^{(n)}.h_1)...).h_n$.

Thus one obtains the differential equation

(23) $(du(x)/dx).h_n = (...(y^{(n)}.h_1)...).h_n = g(u)$.

This equation is with separated variables (see also §4) and over octonions it is equivalent to the equality:

(24) $(1/g(u))(du(x)/dx).h_n = 1$.

Integrating the latter equation we infer:

(25) $\int (1/g(u))du(\omega_{h_n}(x)) = x + \phi_1(x^\circ)$.

This integral becomes simpler, when there exists a super-differentiable diffeomorphism $\nu$ of the domain $U$ so that $\nu(h_n(x))$ is real for each $x$. Suppose that after such integration $u$ is expressed in the form: $u = f(x, \phi_1(x^\circ))$. After $(n-1)$-times non-commutative line integration one gets $y(x)$ (see Equations (3, 4) above).

When the differential equation is:

(26) $(...(y^{(n)}(x).h_1)...).h_n = g((...(y^{(n-2)}(x).h_1)...).h_{n-2})$

with the substitution

(27) $u = (...(y^{(n-2)}.h_1)...).h_{n-2}$ reduces it to the second order differential equation. Then

(28) $(du(x)/dx).h_{n-1} = (...(y^{(n-1)}(x).h_1)...).h_{n-1} =: v(x)$ and

(29) $(dv(x)/dx).h_n = g(u)$.

Now we multiply both sides of the latter equation on $v$ from the left:

(30) $v[(dv(x)/dx).h_n] = [(du(x)/dx).h_{n-1}]g(u)$ and from the right:

(31) $[(dv(x)/dx).h_n]v = g(u)[(du(x)/dx).h_{n-1}]$.

On the other hand, the (super-)derivative of $v^2$ is the following: $(dv^2/dx).h_n = v[(dv(x)/dx).h_n] + [(dv(x)/dx).h_n]v$. Therefore, together with (30, 31) this implies:

(32) $(dv^2/dx).h_n = [(du(x)/dx).h_{n-1}]g(u) + g(u)[(du(x)/dx).h_{n-1}]$.

Integrating Equation (32) over the Cayley-Dickson algebra $\mathcal{A}_r$ we deduce that:

(33) $v^2(\omega_{h_n}(x)) = \int [du(\omega_{h_{n-1}}(x))]g(u) + \int g(u)[du(\omega_{h_{n-1}}(x))] - \phi_{v^2}(Im\,\omega_{h_n})$.

Taking the square root over the Cayley-Dickson algebra from both sides of (33) one obtains:

(34) $v(\omega_{h_n}(x)) = \{\int [du(\omega_{h_{n-1}}(x))]g(u) + \int g(u)[du(\omega_{h_{n-1}}(x))] - \phi_{v^2}(Im\ \omega_{h_n})\}^{1/2}$.

The separation of variables leads to:

(35) $x = \int \{\int [du(\omega_{h_{n-1}}(x))]g(u) + \int g(u)[du(\omega_{h_{n-1}}(x))] - \phi_{v^2}(Im\ \omega_{h_n})\}^{-1/2} du(\omega_{h_n}(x)) - \phi_1(x^\diamond)$,

since $v(\omega_{h_n}(x)) = (du(\omega_{h_n}(x))/dx).1 = (du(\omega_{h_n})/d\omega_{h_n}).h_n(x)$, where $\omega_{h_n} = \omega_{h_n}(x)$. If $u$ is expressed from the latter equality, one gets the general solution $y$ due to Formulas $(3,6)$ above applied to $(...(y^{(n-2)}.h_1)...).h_{n-2}$ instead of $(...(y^{(n)}.h_1)...).h_n$.

**14. Note.** It is important to mention simplifications for resolving differential equations over Cayley-Dickson algebras

(1) $\sum_j a_j(x)y^{(1)}(x).h_j(x) = y^{(1)}(x).[\sum_j a_j(x)h_j(x)]$,

when $a_j(x) \in \mathbf{R}$ and $h_j(x) \in \mathcal{A}_r$ for each $x \in U$, since the real field is the center $Z(\mathcal{A}_r)$ of the Cayley-Dickson algebra $\mathcal{A}_r$, $2 \le r$. This means that the system of partial differential equation with generally non-constant coefficients $a_j$:

(2) $\sum_j a_j(x)y^{(1)}(x).h_j(x) + g(x,y) = 0$

reduces to the ordinary differential equation over the Cayley-Dickson algebra $\mathcal{A}_r$:

(3) $y^{(1)}(x).h(x) + g(x,y) = 0$,

where $h(x) := \sum_j a_j(x)h_j(x)$.

For integration of differential equations over the Cayley-Dickson algebra $\mathcal{A}_r$ the additivity of the non-commutative line integral is helpful:

(4) $\int_\gamma \sum_{j=1}^n f_j(x)dx = \sum_{j=1}^n \int_\gamma f_j(x)$, $n \in \mathbf{N}$.

Evidently one also has:

(5) $\sum_{s;j_1,...,j_n} y^{(n)}.(h_{1,j_1},...,h_{n,j_n})a_{s;1,j_1}(x)...a_{s;n,j_n}(x) = \sum_s y^{(n)}(x).(h_1^s,...,h_n^s)$,

where $h_k^s(x) = \sum_j a_{s;k,j}(x)h_{k,j}(x)$, when $a_{s;j,k}(x) \in \mathbf{R}$ and $h_{j,k}^s(x) \in \mathcal{A}_r$ for all $s,j,k$ and $x \in U$. The sum on the right of (5) may be simpler integrated with the help of (4). When $s$ takes only one value the sum by $s$ certainly gives only one term $y^{(n)}.(h_1,...,h_n)$.

If $A_n$ is a differential operator of the $n$-th order it may happen that there exists the first order differential operators $\kappa_1,...,\kappa_n$ such that $(A_n -$ $\kappa_n...\kappa_1$) is the differential operator $A_{n-1}$ of order not greater than $(n-1)$. Sometimes such reduction may serve for a simplification of integration of the corresponding differential equation

(6) $A_n y = g$ putting

(7) $\kappa_n...\kappa_1 y = f$ and

(8) $A_{n-1} y = g - f$.

Equation (7) can be frequently integrated as the sequence of differential equations of the first order using their types presented above.

Mention, that differential operators of the form $(\nabla u, w)$ with real-valued $u$, $\nabla u = (\partial u/\partial z_1, ..., \partial u/\partial z_m)$ and real-valued (super-)differentiable functions $w = w_1(z),...,w_m(z)$ is of the form given by Equation (1), where $(\nabla u, w) := \sum_{j=1}^{m} w_j(z) \partial u(z)/\partial z_j$. Therefore, equations of the form

(9) $\partial^k u(z)/\partial z_0^k + (\nabla u(z), w(z)) = g(z)$

with $k \in \mathbf{N}$ reduce to that of considered above.

**15. Example.** Let us consider the differential equation

(1) $y^{(2)}(x).(h_1, h_2) + [(dy(x)/dx).h_1]f(x) + [(dy(x)/dx).h_1]g(y)[(dy(x)/dx).h_1] = 0$ over quaternions for super-differentiable functions $h_1 = h_1(x)$, $h_2 = const \in \mathbf{H}$, $f(x)$ and $g(y)$.

To integrate this differential equation we put

(2) $[(dy(x)/dx).h_1] = v(y)h_1 u(x)$. Therefore,

(3) $y^{(2)}(x).(h_1, h_2) = v(y)h_1[(du(x)/dx).h_2] + \{(dv(y)/dy).[(dy(x)/dx).h_2]\}h_1 u = (v(y)h_1)[(du/dx).h_2] + \{(dv(y)/dy).(v(y)h_2 u(x))\}h_1 u(x)$.

And the differential equation takes the view:

(4) $v(y)h_1[(du(x)/dx).h_2] + \{(dv(y)/dy).(v(y)h_2 u(x))\}(h_1 u(x)) + v(y)h_1 u(x)f(x) + (v(y)h_1 u(x))g(y)(v(y)h_1 u(x)) = 0$.

Suppose that either $u(x) \in \mathbf{R}$ or $h_1 \in \mathbf{R}$ or $v(y) \in \mathbf{R}$. The second demand can be fulfilled with a change of variables if use $\omega_{h_1}(x)$. The first or the third condition means that $y^{(1)}$ is either right or left super-linear. Thus the second demand may happen to be the less restrictive between these three conditions. Take the case $u \in \mathbf{R}$ and $h_1 = 1 = i_0$, consequently,

(5) $v[\{(du/dx).h_2 + uf\}/u^2 + v[(1/v)\{(dv/dy).v\}h_2 + gv] = 0$.

Equation (1) certainly has the trivial solution $y' = 0$, when $u(x) = 0$ for each $x$. To obtain a non-trivial solution we put

(6) $(du/dx).h_2 + uf = 0$ and

(7) $(1/v)\{(dv/dy).(vh_2)\} + gv = 0$.

Equation (6) leads to:

(8) $\int (1/u(x))[(du(x)/dx).h_2]dx = -\int f(x)dx$ or

(9) $u(\omega_{h_2}(x)) = \phi_u(Im\ \omega_{h_2})\exp(-\int f(x)dx)$.

The condition $u \in \mathbf{R}$ can be fulfilled for $f \in \mathbf{R}$ with a suitable function $\phi_u(Im\ \omega_{h_2})$. At the same time the integration of Equality (7) gives:

(10) $\int (1/v)[(dv/dy).(v(y)h_2)](1/v)dy = -\int g(y)dy$.

From the identity $[d(v/v)/dy].h = 0$ and the Leibnitz rule over Cayley-Dickson algebras one gets $[d(1/v)/dy].h = -(1/v)[(dv(y)/dy).h](1/v)$. Thus the right side of (10) can be integrated and this implies:

(11) $v(\omega_{vh_2}(y)) = -1/[\phi_{1/v}(Im(\omega_{vh_2}(y))) + \int g(y)dy]$.

Generally relation (11) is rather complicated. For $h_2 = 1$ and $v \in \mathbf{R}$ Equation (11) simplifies to $v = C\exp(-\int g(y)dy)$, since in such variant $(1/v)[(dv/dy).(v(y)h_2)](1/v) = (1/v)(dv/dy)$, $C = C(Im(y))$, $Im(y) := y - Re(y)$. The demand $v \in \mathbf{R}$ can be met with the real-valued mapping $g$ with the corresponding function $C(Im\ y)$.

**16. Example. Transition operator.** In the partial differential equation

(1) $\partial \mathcal{E}/\partial t + (s, \nabla \mathcal{E}) + \alpha \mathcal{E} = \delta(t, x)$

the operator on the right side is called the transition operator, where $0 \leq t \in \mathbf{R}$, $x \in \mathbf{R}^3$, $\alpha$ is a real constant, $\mathcal{E}$ is a generalized function. We suppose that $s$ is a constant vector of unit norm $|s| = 1$, $z = z_0 i_0 + ... + z_3 i_3 \in \mathbf{H}$, $z_0 = t$, $z_j = x_j$ for each $j = 1, 2, 3$, where $\nabla f := (\partial f/\partial x_1, \partial f/\partial x_2, \partial f/\partial x_3)$ denotes the gradient of the differentiable function $f$,

$(x, y) = x_1 y_1 + x_2 y_2 + x_3 y_3$ denotes the scalar product in the Euclidean space $\mathbf{R}^3$. Let us seek a solution in the form:

(2) $\mathcal{E}(t, x) = \theta(t)e^{-\alpha t}u(t, x)$,

where $\theta$ and $u$ are generalized functions. Then Equation (1) simplifies to

(3) $\theta^{(1)}e^{-\alpha t}u + \theta e^{-\alpha t}\partial u/\partial t + \theta e^{-\alpha t}(s, \nabla u) = \delta$ or

(4) $\theta^{(1)}e^{-\alpha t} = -\theta e^{-\alpha t}[(\partial u/\partial t)/u + (s, \nabla u)/u] + \delta/u$,

when $u = u(t, x) \neq 0$. The left side of (4) is independent of $x$, consequently, the right side also is independent of $x$. This means that $[-(\partial u/\partial t) - (s, \nabla u) +$ $\delta/(\theta e^{-\alpha t})]$ is independent of $x$. We put $h = i_0 + s_1 i_1 + s_2 i_2 + s_3 i_3$ and consider $u$ as the function $u = u(z)$ of the quaternion variable $z \in \mathbf{H}$. Therefore,

$$\int_0^x [-(\partial u/\partial t)-(s,\nabla u)+\delta e^{\alpha t}/\theta]dz_0 = -u(\omega_h(x))-\phi_u(Im\ \omega_h(x))-\int [\delta(z)e^{\alpha t}/\theta]dz_0.$$

The equality $\omega_h(y) = x$ gives the solution $y_j = x_j - ts_j$ for each $j = 1, 2, 3$, consequently, $u(t, x) = u(x - ts)$.

One can take any sequence of $z$-differentiable functions $g_n$ with compact supports tending to $\delta(z)$, when $n$ tends to the infinity. That is

$$\lim_{n\to\infty} \int_{-\infty}^{\infty} g_n(t)s(t)dt = s(0) = (\delta, s)$$

for any continuous function $s$ on $\mathbf{R}$. Substituting this expression for $u$ into (3) on gets the simplified differential equation for $z \neq 0$:

(5) $\theta^{(1)}e^{-\alpha t} = -\theta e^{-\alpha t}[-su^{(1)} + su^{(1)}] = 0.$

It is seen that the latter equation has the solution $\theta(t) = 1$ for $t > 0$ and $\theta(t) = 0$ for $t < 0$. At the same time this function $\theta$ as the generalized function has the derivative $\theta^{(1)}(t) = \delta(t)$, consequently, (3) and (5) imply $u(x) = \delta(x)$ and inevitably $u = \delta(x - ts)$. Thus

(6) $\mathcal{E}(t, x) = \theta(t)e^{-\alpha t}\delta(x - ts)$

and the generalized solution of the transition equation

$\partial\mathcal{E}/\partial t + (s, \nabla\mathcal{E}) + \alpha\mathcal{E} = p$ is the following $\mathcal{E} * p$, where $p = p(z)$ is a generalized function or particularly a $z$-differentiable function, $z \in \mathbf{H}$.

**20. Remark.** The results presented above show that the utilization of the Cayley-Dickson algebras enlarges possibilities for the integration of differential equations and their systems in comparison with the real and complex cases. With the help of generalized functions the integration methods of differential equations with octonion variables can be spread on Sobolev spaces. Dirac had used the complexified bi-quaternion algebra to solve Klein-Gordon equation for quantum particles with spin. This partial differential equation is of the hyperbolic type with constant coefficients. Using the technique of this article it is possible to integrate subsequently hyperbolic partial differential equations with variable coefficients (see also §§13 and 14 above).

This is planned to be written in details in subsequent papers together with differential equations of generalized functions of octonion variables, singular integrals (exceptional solutions) and examples of differential equations and their fundamental and general solutions.

Department of Applied Mathematics,

Moscow State Technical University MIREA, av. Vernadsky 78,

Moscow, Russia

e-mail: sludkowski@mail.ru